\documentclass[leqno,final]{siamltex}
\usepackage{graphicx} 
\usepackage{amsmath,amstext,amssymb,bm}
\usepackage{leftidx}

\usepackage{xcolor} 
\usepackage{soul} 
\usepackage{tikz}
\usetikzlibrary{shapes,arrows}
\usepackage{mathrsfs}
\usepackage{epstopdf}
\usepackage{color}
\usepackage{multirow}
\usepackage{tabularx}
\usepackage[shortlabels]{enumitem}

\numberwithin{equation}{section}

\allowdisplaybreaks[4]

\newcommand{\cP}{\mathcal{P}}

\newcommand{\cF}{\mathcal{F}}

\newcommand{\mP}{\mathbb{P}}

\newcommand{\mE}{\mathbb{E}}

\newcommand{\Ome}{\Omega}
\newcommand{\p}{\partial}
\newcommand{\nab}{\nabla}

\newcommand{\vb}{{\bf b}}

\newcommand{\vH}{{\bf H}}

\newcommand{\vv}{{\bf v}}

\newcommand{\e}{\pmb{\varepsilon}}

\newcommand{\pphi}{\pmb{\phi}}

\newcommand{\vvarphi}{\pmb{\varphi}}
\newcommand{\vsigma}{\pmb{\sigma}}
\newcommand{\vtheta}{\pmb{\theta}}

\begin{document}
	
	\title{Analysis of fully discrete Crank-Nicolson finite element methods for a stochastic Keller-Segel chemotaxis system with gradient type multiplicative noise\thanks{This work was partially supported by the NSF grant DMS-2530211.}}
	\markboth{LIET VO}{STOCHASTIC KELLER-SEGEL CHEMOTAXIS EQUATIONS}		
	
	\author{
		Liet Vo\thanks{School of Mathematical and Statistical Sciences, The University of Texas Rio Grande Valley, Edinburg, TX 78539, U.S.A.  ({\tt liet.vo@utrgv.edu}).} 
	}

	\maketitle
	
	\begin{abstract}
		We develop and analyze numerical methods for a stochastic Keller–Segel system perturbed by Stratonovich noise, which models chemotactic behavior under randomly fluctuating environmental conditions. The proposed fully discrete scheme couples a Crank–Nicolson time discretization with a splitting mixed finite element method in space. We rigorously prove the stability of the numerical scheme and establish strong convergence rates of order $O(k^{1/2} + k^{-1/2}h^2)$, where $k$ and $h$ denote the time and spatial step sizes, respectively. Notably, the presence of stochastic forcing leads to an inverse dependence on $k$ in the error estimates, distinguishing the convergence behavior from that of the deterministic case. Numerical experiments are presented to validate the theoretical results and demonstrate the effectiveness and accuracy of the proposed methods.

	\end{abstract}
	
	
	
	

	\section{Introduction}\label{sec-1}
	In this paper, we consider the following stochastic Keller-Segel system with gradient-type multiplicative noise:
	\begin{subequations}\label{eq1.1}
		\begin{align} \label{eq1.1a}
			du  &=  \bigl[\nu\Delta u -\chi\nab\cdot (u \nab c)  \bigr]\, dt + \delta\, {\bf b}\cdot\nab u\circ dW(t) \quad\mbox{a.s. in}\, D_T:=(0,T)\times D,\\ 
			\label{eq1.1b}	\Delta c&=  c-u\qquad\mbox{a.s. in}\, D_T,    \\
			u(0) &= u_0, \quad c(0) = c_0 \qquad\mbox{a.s. in}\, D.
		\end{align}
	\end{subequations}
	where $D=[0,L]^2\subset \mathbb{R}^2$, $L>0, T > 0$. $u$ and $c$ denote, respectively, the cell density and the concentration of an attractive chemical signal, which are periodic in space with period $L$. $\nu, \delta$ are the positive constants which represent the cell diffusion constant and the intensity of the noise. $\vb: \mathbb{R}^2 \rightarrow \mathbb{R}^2$ is a vector constant. The given constant $\chi>0$ denotes the chemotactic sensitivity. $\{W(t): t\geq 0\}$ is a real Wiener process. The notation $``\circ"$ denotes the Stratonovich noise. 
	
	Chemotaxis refers to the directed movement of cells or organisms in response to chemical gradients in their environment. It is a fundamental mechanism in numerous biological processes, including immune responses, wound healing, embryonic development, and tumor invasion~\cite{friedl2009collective, roussos2011chemotaxis}. A well-known example is the movement of Escherichia coli bacteria toward higher concentrations of nutrients, or away from harmful substances~\cite{berg1972chemotaxis}. At the cellular level, chemotactic movement is regulated by signaling pathways that enable cells to detect and respond to chemoattractant or chemorepellent cues~\cite{devreotes2003chemotaxis}.
	
	Mathematically, chemotaxis is often modeled using systems of partial differential equations (PDEs), with the Keller–Segel model \cite{keller1970initiation} (the equations \eqref{eq1.1} with $\delta =0$) being one of the most widely studied frameworks. This model describes the spatiotemporal evolution of cell density and chemical concentration, accounting for diffusion and chemotactic drift. Chemotaxis models exhibit complex behaviors such as aggregation, blow-up, and pattern formation~\cite{horstmann2003limits, perthame2006transport}. Understanding the analytical and numerical properties of these systems is essential for capturing the underlying biological phenomena and for designing reliable computational tools to simulate them. 
	
	Numerical methods for deterministic Keller–Segel systems remain an active area of research, with a variety of approaches such as finite element methods~\cite{shi2023linearized, saito2012error, chen2022error}, finite volume methods~\cite{filbet2006finite, zhou2017finite}, and finite difference methods~\cite{hu2023positivity}. A central challenge in designing numerical methods for these systems lies in accurately capturing key features of the true solution, such as finite-time blow-up, positivity preservation, mass conservation, and energy decay in the approximate solution~\cite{nagai2000behavior}. Additionally, the nonlinear term in~\eqref{eq1.1a}, which is non-Lipschitz, poses significant difficulties in the context of error analysis.
	
	Recently, stochastic Keller–Segel systems have drawn increasing attention due to their relevance in modeling real-world phenomena. In the derivation of deterministic models, fluctuations around the mean behavior are often neglected. However, more realistic models must account for inherent randomness and variability in the natural environment, which are not reproducible~\cite{hausenblas2022one}. Consequently, incorporating noise into the system is essential for more accurate modeling.
	
	The well-posedness of the stochastic Keller–Segel system~\eqref{eq1.1} has been studied in~\cite{huang2021microscopic, misiats2022global}. Notably, these works establish the existence and uniqueness of pathwise solutions~\cite[Theorem 3.3]{huang2021microscopic}. Furthermore, it has been shown that the solution $u$ can blow up in finite time when the initial data is sufficiently large~\cite[Theorem 4.2]{misiats2022global}, specifically when
	\begin{align}
		\chi \int_D u_0(x)\, dx \geq 4\pi(2\nu + \delta^2).
	\end{align}
	
	In addition, the density solution $u$ preserves positivity and mass conservation, although energy decay has not yet been established~\cite[Proposition 2.4]{misiats2022global}. While there has been progress on the theoretical PDE of the system~\eqref{eq1.1}, the development of numerical methods with rigorous analysis remains largely open. The primary challenges include: (i) ensuring the numerical method faithfully captures key properties of the solution, and (ii) handling the interplay between the nonlinear second-order term and the stochastic forcing.
	
	To the best of our knowledge, no numerical methods have yet been established for solving~\eqref{eq1.1}. In this paper, we address this gap by proposing and analyzing a Crank-Nicolson method for time discretization combined with a splitting mixed finite element method for spatial discretization. The Crank-Nicolson method is particularly well suited for approximating Stratonovich-type noise. On the other hand, the nonlinear term in~\eqref{eq1.1a} is both non-Lipschitz and sign-indefinite, creating significant challenges for traditional stability and error analysis. To overcome this, we adopt a splitting mixed finite element approach inspired by~\cite{zhang2016characteristic, duarte2021numerical}, where we introduce an auxiliary variable $\vsigma = \nabla c$ to control the second-order nonlinear term in the density equation.
	
	An important observation is that the stochastic nature of the problem influences the error estimates of the numerical solution. As a result, the error bounds, $\mathcal{O}(k^{1/2} + h + k^{-1/2}h^2)$, include a factor of $k^{-1/2}$, which does not appear in deterministic counterparts. We rigorously prove that these error estimates are sharp and confirm them through numerical experiments. We also verify that our numerical solution satisfies the mass conservation. Additionally, we demonstrate numerically that the approximate density remains nonnegative and exhibits finite-time blow-up when the initial data is sufficiently large.
	
	The remainder of the paper is organized as follows. In Section~\ref{sec2}, we introduce the necessary notation, function spaces, and preliminary results. We also present the regularity of the solution and derive the mixed formulation of~\eqref{eq1.1}. Section~\ref{sec3} is devoted to the Crank-Nicolson semi-discrete scheme, where we establish stability properties and derive strong error estimates. In Section~\ref{sec4}, we present the fully discrete splitting mixed finite element method and analyze its error bounds. Finally, in Section~\ref{sec5}, we provide a series of numerical experiments to validate the theoretical findings, and Section \ref{sec6} gives a conclusion of the paper.

	\section{Prelimimaries}\label{sec2}
	\subsection{Notations and assumptions}
	Standard function and space notation will be adopted in this paper. 
	Let $H^1_0(D)$ denote the subspace of $H^1(D)$ whose ${\mathbb R}^d$-valued functions have zero trace on $\p D$, and $(\cdot,\cdot):=(\cdot,\cdot)_D$ denote the standard $L^2$-inner product, with induced norm $\|\cdot\|_{L^2}$. We also denote ${ L}^p_{per}(D)$ and ${H}^{k}_{per}(D)$ as the Lebesgue and Sobolev spaces of the functions that are periodic with period $L$ in each coordinate direction for almost every ${x} \in D$ { and have vanishing mean.} In addition, we denote $\vH_{per}^{k}(D) = (H_{per}^k(D), H_{per}^k(D))$.
	$C$ denotes a generic constant that is independent of the mesh parameters $h$ and $k$.
	
	Let $(\Omega,\cF, \{\cF_t\},\mP)$ be a filtered probability space with the probability measure $\mP$, the 
	$\sigma$-algebra $\cF$ and the continuous  filtration $\{\cF_t\} \subset \cF$. For a random variable $v$ 
	defined on $(\Omega,\cF, \{\cF_t\},\mP)$,
	${\mathbb E}[v]$ denotes the expected value of $v$. 
	For a vector space $X$ with norm $\|\cdot\|_{X}$,  and $1 \leq p < \infty$, we define the Bochner space
	$\bigl(L^p(\Omega;X); \|v\|_{L^p(\Omega;X)} \bigr)$, where
	$\|v\|_{L^p(\Omega;X)}:=\bigl({\mathbb E} [ \Vert v \Vert_X^p]\bigr)^{\frac1p}$.
	
	We also recall the usefull elliptic estimate $\|u\|_{H^2} \leq C(D)\|\Delta u\|_{L^2}$, where $C(D)>0$ and only depends on the domain $D$, from  \cite[Theorem 9.3]{gilbarg1977elliptic}.
	
	Moreover, the Ladyzhenskaya inequality in $\mathbb{R}^2$ that is $\|u\|_{L^4}\leq C_L\|u\|^{1/2}_{L^2}\|\nab u\|^{1/2}_{L^2}$, will be used in the paper.
	
	Lastly, we use the equivalent norms (see \cite[Corollary 3.5]{amrouche2013lp} and \cite[Section 2.1]{duarte2021numerical}):
	\begin{align}
		\|\vsigma\|^2_{H^1} = \|\vsigma\|^2_{L^2} + \|\nab\cdot \vsigma\|^2_{L^2} + \|rot \vsigma\|^2_{L^2}\qquad\forall\vsigma \in \vH^1_{per}(D).
	\end{align}
	\subsection{Useful results}
	We state a useful nonlinear Gronwall lemma, which will be used later to establish stability estimates of the proposed methods. 
	\begin{lemma}\cite[Theorem 2.3.5]{pachpatte2001inequalities}\label{nonlinear_Gronwall}
		Let $f(m)$, $b(m)$ and $k(m)$ be nonnegative functions defined for $m \in \mathbb{N}\cup\{0\}$ and
		\begin{align*}
			f(m) \leq a + \sum_{n=0}^{m-1} b(n)f(n) + \sum_{n=0}^{m-1} k(n) f^p(n),
		\end{align*}
		where $a>0$ and $p\geq 0$, $p\neq 1$, are constants. Then,
		\begin{align*}
			f(m) \leq \zeta^{-1}(m)\left\{a^{1-p} + (1-p)\sum_{n=0}^{m-1} k(m)\zeta^{1-p}(n+1)\right\}^{1/(1-p)},
		\end{align*}
		for $m \in \mathbb{N}$ such that 
		\begin{align*}
			a^{1-p} + (1-p)\sum_{n=0}^{m-1}k(s)\zeta^{1-p}(n+1) >0, 
		\end{align*}
		and $\zeta(m) = \prod_{n=0}^{m-1} (1 + b(n))^{-1}$.
	\end{lemma}

	\subsection{Solution concepts}
	In \cite{mayorcas2023blow,misiats2022global}, the author established a weak solution for the system \eqref{eq1.1} that satisfies the following identities:
	\begin{subequations}\label{eq_2.1}
		\begin{align}\label{eq_2.1a}
			(u(t),\phi) &= (u_0,\phi) - \nu \int_0^t (\nab u(s), \nab \phi)\, ds \\\nonumber
			&\qquad+ \chi \int_0^t (u\nab c, \nab \phi)\, ds + \delta \int_0^t (\vb\cdot\nab u(s),\phi)\circ dW(s),\\
			\label{eq_2.1b}		(\nab c(t),& \nab\psi) + (c(t), \psi) = (u(t),\psi),
		\end{align}
	\end{subequations}
	for all $\phi \in H_{per}^1(D)$, and $\psi \in H_{per}^1(D)$.
	
	In addition, using the conversion between the It\^o integrals and the Stratonovich integrals, we can rewrite the system \eqref{eq2.1} as follows.
	\begin{subequations}\label{eq_2.2}
		\begin{align}\label{eq_2.2a}
			(u(t),\phi) &= (u_0,\phi) - \nu \int_0^t (\nab u(s), \nab \phi)\, ds + \chi \int_0^t (u\nab c, \nab \phi)\, ds \\\nonumber
			&\qquad+\frac{\delta^2}{2}\int_0^t(\nab\cdot(B\nab u(s)), \phi)\, ds + \delta \int_0^t (\vb\cdot\nab u(s),\phi)\, dW(s),\\
			\label{eq_2.2b}	(\nab c(t),& \nab\psi) + (c(t), \psi) = (u(t),\psi),
		\end{align}
	\end{subequations}
	where $B = \vb\otimes\vb\in \mathbb{R}^{2\times2}$ and $B_{ij} = b_ib_j$ for all $i, j = 1,2$.
	
	When $u_0$ is smooth enough, there exists a unique solution \cite[Theorem 3.3]{huang2021microscopic} of \eqref{eq_2.2} that satisfies the following regularity:
	\begin{lemma}\label{Stability_PDE} Let $(u,c)$ be the solution of \eqref{eq_2.1}. Then, the following estimates satisfy:
		\begin{enumerate}
			\item[(a)] For $p\geq 2$, $\displaystyle \|u\|_{L^p} \leq \|u_0\|_{L^p}$ $\mP$-a.s.
			\item[(b)] Suppose that for $s\geq 0$, if $u \in H_{per}^s(D)$, then there exists a constant $C>0$ such that
			\begin{align*}
				\|c(t)\|_{H^{s+2}} \leq C\|u(t)\|_{H^s} \quad\mP-a.s., \quad\forall t \in [0,T].
			\end{align*}
			\item[(c)] If $u_0 \in H_{per}^1(D)$ and suppose that $\exists \delta_0>0$ such that $\frac{7\nu}{4} - 2\delta^2C(D)|\vb|^2>0$ for all $\delta \in (0,\delta_0)$, then there exists some constant $C_1 \equiv C(\chi,u_0,\nu, T)>0$ such that 
			\begin{align*}
				\sup_{t\in [0,T]}\|\nab u(t)\|^2_{L^2} + 
				\left(\frac{7\nu}{4} - 2\delta^2C(D)|\vb|^2\right)\int_0^T\|\Delta u(s)\|^2_{L^2}\, ds \leq C_1, \qquad\mP-a.s..
			\end{align*}
			\item[(d)] If $u_0 \in H^2_{per}(D)$ and suppose that $\exists \delta_0>0$ such that $\frac{7\nu}{4} - 2\delta^2C(D)|\vb|^2>0$ for all $\delta \in (0,\delta_0)$, then there exists some constant $C_2 \equiv C(\chi,u_0,\nu, T, C_1)>0$ such that 
			\begin{align*}
				\sup_{t\in [0,T]}\|\Delta u(t)\|^2_{L^2} + 
				\left(\frac{7\nu}{4} - 2\delta^2C(D)|\vb|^2\right)\int_0^T\|\nab\Delta u(s)\|^2_{L^2}\, ds &\leq C_2,\qquad\mP-a.s..
			\end{align*}
		\end{enumerate}
	\end{lemma}
	\begin{proof} First, Part (a) has been proved in \cite[Proposition 2.5, Theorem 3.2]{misiats2022global}. Part (b) is the standard elliptic theory, which can be found easily in many PDE textbooks. Proving Parts (c) and (d) is very similar, so we focus on establishing Part (c) and leave Part (d) as an exercise for the reader.
		
		Applying the It\^o formula to $\Phi(u(t)) = \|\nab u(t)\|^2_{L^2}$, we have
		\begin{align}\label{eq_2.3}
				\|\nab u(t)\|^2_{L^2} + 2\nu \int_{0}^t \|\Delta u(\xi)\|^2_{L^2}\, d\xi
			&=   
			-2\chi\int_{0}^t (\nab u(\xi),\nab(\nab\cdot(u(\xi)\nab c(\xi))))\, d\xi\\\nonumber
			&\qquad + 2\delta \int_{0}^t (\nab u(\xi),\nab(\vb\cdot\nab u(\xi)))\, dW(\xi)\\\nonumber
			&\qquad+\delta^2 \int_{0}^t (\nab u(\xi),\nab(\nab\cdot(B\nab u(\xi))))\, d\xi\\\nonumber
			&\qquad+ \delta^2 \int_0^t \|\nab(\vb\cdot\nab u(\xi))\|^2_{L^2}\, d\xi + \|\nab u_0\|^2_{L^2}\\\nonumber
			&:=L_1 + \cdots + L_4 + \|\nab u_0\|^2_{L^2}.
		\end{align}
		
		First, we notice that by the divergence theorem,
		\begin{align*}
			L_2 = 2\delta \int_{0}^t (\nab u(\xi),\nab(\vb\cdot\nab u(\xi)))\, dW(\xi) = 0.
		\end{align*}
		
		Using integration by parts, we get
		\begin{align*}
			L_1 &= 2\chi\int_{0}^t (\Delta u(\xi), \nab\cdot(u(\xi)\nab c(\xi)))\, d\xi\\\nonumber
			&= 2\chi\int_{0}^t (\Delta u(\xi), \nab u(\xi)\nab c(\xi) + u \Delta c(\xi))\, d\xi\\\nonumber
			&\leq \frac{\nu}{4} \int_{0}^t \|\Delta u(\xi)\|^2_{L^2}\, d\xi + \frac{C\chi^2}{\nu} \int_{0}^t (\|\nab u(\xi)\nab c(\xi)\|^2_{L^2} + \|u\Delta c(\xi)\|^2_{L^2})\, d\xi\\\nonumber
			&\leq \frac{\nu}{4} \int_{0}^t \|\Delta u(\xi)\|^2_{L^2}\, d\xi + \frac{C\chi^2}{\nu} \int_{0}^t (\|\nab u(\xi)\|^2_{L^4}\|\nab c(\xi)\|^2_{L^4} + \|u\|^2_{L^{4}}\|\Delta c(\xi)\|^2_{L^4})\, d\xi.
		\end{align*}
		
		Next, using the Gagliardo–Nirenberg inequality, Lemma \ref{Stability_PDE} Part (a) and Part (b), we have
		\begin{align*}
			&L_1 -\frac{\nu}{4} \int_{0}^t \|\Delta u(\xi)\|^2_{L^2}\, d\xi\\\nonumber
			&\leq  \frac{C\chi^2}{\nu} \int_{0}^t (\|\nab u(\xi)\|_{L^2}\|\Delta u(\xi)\|_{L^2}\| c(\xi)\|^2_{H^2} + \|u\|_{L^{2}}\|\nab u(\xi)\|_{L^2}\|\Delta c(\xi)\|_{L^2}\|c(\xi)\|_{H^3})\, d\xi\\\nonumber
			&\leq \frac{C\chi^2}{\nu} \int_{0}^t (\|\nab u(\xi)\|_{L^2}\|\Delta u(\xi)\|_{L^2}\|u_0\|^2_{L^2} + \|u_0\|^2_{L^{2}}\|\nab u(\xi)\|_{L^2}^2)\, d\xi\\\nonumber
			&\leq \frac{\nu}{4}\int_{0}^t \|\Delta u(\xi)\|^2_{L^2}\, d\xi +  \frac{C\chi^4}{\nu^3} \int_{0}^t \|u_0\|^2_{L^2}(1+\|u_0\|^2_{L^2})\|\nab u(\xi)\|^2_{L^2}\, d\xi\\\nonumber
		\end{align*}
		In addition, using the inequality $\|u\|_{H^2} \leq C(D)\|\Delta u\|_{L^2}$, we can control $L_3$ and $L_4$ as follows:
		\begin{align*}
			L_3 &= -\delta^2\int_0^t \bigl(\Delta u(\xi), \nab\cdot(B\nab u(\xi))\bigr)\, d\xi\\\nonumber
			&\leq \delta^2\int_0^t \|\Delta u(\xi)\|_{L^2} \|\nab\cdot(B\nab u(\xi))\|_{L^2}\, d\xi\\\nonumber
			&\leq \delta^2\int_0^t \|\Delta u(\xi)\|_{L^2}|\vb|^2 \|u(\xi)\|_{H^2}\, d\xi\\\nonumber
			&\leq \delta^2C(D)|\vb|^2\int_0^t \|\Delta u(\xi)\|^2_{L^2}\, d\xi,
		\end{align*}
		where $C(D)$ depends only on the domain $D$.
		
		Similarly, we also have
		\begin{align*}
			L_4 &\leq \delta^2 C(D)|\vb|^2\int_{0}^t \|\Delta u(\xi)\|^2_{L^2}\, d\xi.
		\end{align*}
		
		Substituting all the estimates from $L_1, \dots, L_4$ into \eqref{eq_2.3} we obtain
		\begin{align*}
			&	\|\nab u(t)\|^2_{L^2} + \left(\frac{7\nu}{4} - 2\delta^2C(D)|\vb|^2\right)\int_{0}^t \|\Delta u(\xi)\|^2_{L^2}\, d\xi \\\nonumber
			&\leq \|\nab u_0\|^2_{L^2} +  \frac{C\chi^4}{\nu^3} \int_{0}^t \|u_0\|^2_{L^2}(1+\|u_0\|^2_{L^2})\|\nab u(\xi)\|^2_{L^2}\, d\xi,
		\end{align*}
		where $\frac{7\nu}{4} - 2\delta^2C(D)|\vb|^2 >0$ by the assumption of $\delta$.
		
		Finally, the proof is complete by using the Gronwall inequality on the above inequality. 
		
	\end{proof}
	
	\subsection{Variational mixed formulation}
	To control the nonlinear term of \eqref{eq_2.1a}, we first introduce a new variable $\vsigma = \nab c$ in \eqref{eq_2.1b}, we obtain
	\begin{align}\label{eq2.1}
		\bigl(\vsigma,\vvarphi\bigr) + \bigl(c,\nab\cdot \vvarphi\bigr) &= 0,\\
		\label{eq2.2}	\bigl(\nab\cdot\vsigma, \psi\bigr) + (u,\psi)&= \bigl(c,\psi\bigr),
	\end{align}
	for $\vvarphi\in \vH_{per}^1(D)$, and $\psi\in H_{per}^1(D)$. Next, choosing $\psi = \nab\cdot\vvarphi$ in $\eqref{eq2.2}$ and replacing the result into $\eqref{eq2.1}$, and also adding a zero term $\bigl(rot\, \vsigma,\, rot \,\vvarphi \bigr)$ (this is zero because of the fact that $rot \,\vsigma = 0$) to $\eqref{eq2.1}$, we obtain the following mixed variational formulation: Find $\bigl(u,\vsigma, c\bigr) \in{H}^1_{per}(D) \times \vH^1_{per}(D) \times L^2_{per}(D)$ such that
	\begin{subequations}\label{Weak_formulation}
		\begin{align}
			\label{eq_sigma}	&\bigl(\vsigma, \vvarphi\bigr) + \bigl(\nab\cdot \vsigma, \nab\cdot\vvarphi\bigr) + \bigl(rot \,\vsigma, rot \,\vvarphi\bigr) = -\bigl(u,\nab\cdot\vvarphi\bigr)\\
			\label{eq_u}	&\bigl(u(t), \varphi\bigr) + \nu \int_0^t \bigl(\nab u, \nab\varphi\bigr)  \, ds = \bigl(u_0, \varphi\bigr)+\chi \int_0^t \bigl(u\vsigma,\nab \varphi\bigr)\, ds \\\nonumber&\qquad\qquad\qquad\qquad\qquad\qquad\qquad+ \label{eq_c}\delta\Bigl(\int_0^t\vb\cdot\nab u\circ dW(s),\varphi\Bigr),\\
			&\bigl(c,\psi\bigr)=	\bigl(\nab\cdot\vsigma, \psi\bigr) + (u,\psi),
		\end{align}
	\end{subequations}
	for all $\bigl(\varphi,\vvarphi, \psi\bigr) \in {H}^1_{per}(D) \times \vH^1_{per}(D) \times L^2_{per}(D)$.

	\begin{lemma}\label{Lemma_Holder} Let $(u,\vsigma,c)$ be the solution of \eqref{Weak_formulation}. Then, there exist constants $C>0$ such that if
		\begin{enumerate}
			\item[(a)] $u_0 \in L^2(\Ome; H^1_{per}(D))\cap L^4(\Ome; L^2_{per}(D))$, then
			\begin{align}\label{eq_2.6}
				\mE\bigl[\|u(t) - u(s)\|^2_{L^2}\bigr] &\leq C|t-s|,\\
				\label{eq_2.7}	\mE\bigl[\|\vsigma(t) - \vsigma(s)\|^2_{H^1}\bigr] &\leq C|t-s|.
			\end{align}
			\item[(b)] $u_0 \in L^2(\Ome; H^2_{per}(D))\cap L^4(\Ome; H^1_{per}(D))$, then
			\begin{align}
				\label{eq_2.8}	\mE\bigl[\|\nab(u(t) - u(s))\|^2_{L^2}\bigr] &\leq C|t-s|.
			\end{align}
			\item[(c)] $u_0 \in  L^{\infty}(\Omega; H^2_{per}(D))$, then
			\begin{align}
				\label{eq_2.9}	\mE\bigl[\|\nab(u(t) - u(s))\|^4_{L^2}\bigr] &\leq C|t-s|^2.
			\end{align}
		\end{enumerate}
	\end{lemma}
	\begin{proof}
		\begin{enumerate}
			\item[(a)] First, we prove the estimate \eqref{eq_2.6}. Applying the It\^o formula to $\Phi(u(t)) = \|u(t) - u(s)\|^2_{L^2}$, we obtain
			\begin{align*}
				\|u(t) - u(s)\|^2_{L^2} &= -2\nu \int_{s}^{t} (\nab(u(\xi) - u(s)), \nab u(\xi))\, d\xi \\\nonumber
				&\qquad+2 \chi\int_{s}^t (\nab(u(\xi) - u(s)), u(\xi)\nab c(\xi))\, d\xi\\\nonumber
				&\qquad - \delta^2\int_{s}^t (\vb\cdot\nab (u(\xi) - u(s)), \vb\cdot\nab u(\xi))\, d\xi \\\nonumber
				&\qquad + 2\delta\int_{s}^t(\vb\cdot\nab u(\xi), u(\xi)  -u(s))\, dW(\xi) + \delta^2\int_{s}^t \|\vb\cdot\nab u(\xi)\|^2_{L^2}\, d\xi,
			\end{align*}
			which implies that
			\begin{align*}
				\mE[\|u(t) - u(s)\|^2_{L^2}] &= -2\nu \int_{s}^{t} \mE[\|\nab(u(\xi) - u(s))\|^2_{L^2}]\, d\xi \\\nonumber
				&\qquad-2\nu \int_{s}^{t} \mE[(\nab(u(\xi) - u(s)), \nab u(s))]\, d\xi \\\nonumber
				&\qquad+2 \chi\int_{s}^t \mE[(\nab(u(\xi) - u(s)), u(\xi)\nab c(\xi))]\, d\xi\\\nonumber
				&\qquad - \delta^2\int_{s}^t \mE[\|\vb\cdot\nab (u(\xi) - u(s))\|^2_{L^2}]\, d\xi + \delta^2\int_{s}^t \mE[\|\vb\cdot\nab u(\xi)\|^2_{L^2}]\, d\xi\\\nonumber
				&\qquad - \delta^2\int_{s}^t \mE[(\vb\cdot\nab (u(\xi) - u(s)), \vb\cdot\nab u(s))]\, d\xi .
			\end{align*}
			
			Using Cauchy-Schwarz's inequality, we get
			\begin{align*}
				&\mE[\|u(t) - u(s)\|^2_{L^2}] + 2\nu \int_{s}^{t} \mE[\|\nab(u(\xi) - u(s))\|^2_{L^2}]\, d\xi \\\nonumber&\qquad\qquad\qquad\qquad\qquad+  \delta^2\int_{s}^t \mE[\|\vb\cdot\nab (u(\xi) - u(s))\|^2_{L^2}]\, d\xi \\\nonumber
				&\leq C\mE[\|\nab u(s)\|^2_{L^2}]|t-s| + \frac{\nu}{2} \int_{s}^t \mE[\|\nab(u(\xi) - u(s))\|^2_{L^2}]\, d\xi\\\nonumber
				&\qquad+C\sup_{\xi\in [0,T]}\mE[\|u(\xi)\nab c(\xi)\|^2_{L^2}]|t-s| + \frac{\nu}{2} \int_{s}^t \mE[\|\nab(u(\xi) - u(s))\|^2_{L^2}]\, d\xi\\\nonumber
				&\qquad +C\mE[\|\vb\cdot\nab u(s)\|^2_{L^2}]|t-s| + \frac{\delta^2}{2} \int_{s}^t \mE[\|\vb\cdot\nab(u(\xi) - u(s))\|^2_{L^2}]\, d\xi\\\nonumber
				&\qquad+ \delta^2 \sup_{\xi\in [0,T]}\mE[\|\vb\cdot\nab u(\xi)\|^2_{L^2}] |t-s|.
			\end{align*}
			
			Therefore, 
			\begin{align*}
				&\mE[\|u(t) - u(s)\|^2_{L^2}] + 2\nu \int_{s}^{t} \mE[\|\nab(u(\xi) - u(s))\|^2_{L^2}]\, d\xi +  \delta^2\int_{s}^t \mE[\|\vb\cdot\nab (u(\xi) - u(s))\|^2_{L^2}]\, d\xi \\\nonumber
				&\leq C\Bigl\{\sup_{s\in [0,T]}\mE[\|\nab u(s)\|^2_{L^2}] + \sup_{\xi\in [0,T]}\mE[\|u(\xi)\nab c(\xi)\|^2_{L^2}] +\sup_{s\in [0,T]}\mE[\|\vb\cdot\nab u(s)\|^2_{L^2}]\\\nonumber
				&\qquad+\sup_{s\in [0,T]}\mE[\|\vb\cdot\nab u(\xi)\|^2_{L^2}] \Bigr\}|t-s|\\\nonumber
				&\leq C\Bigl\{\sup_{s\in [0,T]}\mE[\|\nab u(s)\|^2_{L^2}] + \sup_{\xi\in [0,T]}\mE[\|u(\xi)\|_{L^2}\|\nab u(\xi)\|_{L^2}\| c(\xi)\|^2_{H^2}] \Bigr\}|t-s|\\\nonumber
				&\leq C\Bigl\{\sup_{s\in [0,T]}\mE[\|\nab u(s)\|^2_{L^2}] + \sup_{\xi\in [0,T]}\mE[\|u(\xi)\|^3_{L^2}\|\nab u(\xi)\|_{L^2}] \Bigr\}|t-s|.
			\end{align*}
			
			Next, using \eqref{eq_2.6} and the equation \eqref{eq_sigma}, we easily establish \eqref{eq_2.7}. The proof is complete.
			
			\item[(b)] Applying the It\^o formula to $\Phi(u(t)) = \|\nab(u(t) - u(s))\|^2_{L^2}$, we obtain
			\begin{align*}
				\|\nab(u(t) - u(s))\|^2_{L^2} &= 2\nu \int_{s}^t (\nab(u(\xi) - u(s)), \nab(\Delta u(\xi)))\, d\xi\\\nonumber
				&\qquad -2\chi\int_{s}^t (\nab(u(\xi) - u(s)),\nab(\nab\cdot(u(\xi)\nab c(\xi))))\, d\xi\\\nonumber
				&\qquad + 2\delta \int_{s}^t (\nab(u(\xi) - u(s)),\nab(\vb\cdot\nab u(\xi)))\, dW(\xi)\\\nonumber
				&\qquad+\delta^2 \int_{s}^t (\nab(u(\xi) - u(s)),\nab(\nab\cdot(B\nab u(\xi))))\, d\xi\\\nonumber
				&\qquad+ \delta^2 \int_s^t \|\nab(\vb\cdot\nab u(\xi))\|^2_{L^2}\, d\xi.
			\end{align*}
			
			Taking the expectation and using integration by parts on the above identity, we arrive at
			\begin{align*}
				&\mE	[\|\nab(u(t) - u(s))\|^2_{L^2}] + 2\nu \int_{s}^t \mE[\|\Delta (u(\xi) - u(s))\|^2_{L^2}]\, d\xi\\\nonumber
				&= -2\nu \int_{s}^t\mE[(\Delta(u(\xi) - u(s)), \Delta u(s))]\, d\xi  \\\nonumber
				&\qquad+2\chi\int_{s}^t \mE[(\Delta(u(\xi) - u(s)), \nab\cdot(u(\xi)\nab c(\xi)))]\, d\xi\\\nonumber
				&\qquad-\delta^2 \int_{s}^t \mE[(\Delta(u(\xi) - u(s)),\nab\cdot(B\nab u(\xi)))]\, d\xi + \delta^2 \int_s^t \mE[\|\nab(\vb\cdot\nab u(\xi))\|^2_{L^2}]\, d\xi.
			\end{align*}
			
			Using Cauchy-Schwarz's inequality, we obtain
			\begin{align*}
				&\mE	[\|\nab(u(t) - u(s))\|^2_{L^2}] + \frac{5\nu}{4} \int_{s}^t \mE[\|\Delta (u(\xi) - u(s))\|^2_{L^2}]\, d\xi\\\nonumber
				&\leq \frac{\nu}{4} \int_{s}^t\mE[\|\Delta(u(\xi) - u(s))\|^2_{L^2}]\, d\xi + C\mE[\|\Delta u(s)\|^2_{L^2}]|t-s| \\\nonumber
				&\qquad+\frac{\nu}{4}\int_{s}^t \mE[\|\Delta(u(\xi) - u(s))\|^2_{L^2}]\, d\xi + C\int_{s}^t\mE[\|\nab\cdot(u(\xi)\nab c(\xi))\|^2_{L^2}]\, d\xi\\\nonumber
				&\qquad +\frac{\nu}{4} \int_{s}^t \mE[\|\Delta(u(\xi) - u(s))\|^2_{L^2}]\, d\xi + C\int_{s}^t\mE[\|\nab\cdot(B\nab u(\xi))\|^2_{L^2}]\, d\xi\\\nonumber
				&\qquad+ \delta^2 \int_s^t \mE[\|\nab(\vb\cdot\nab u(\xi))\|^2_{L^2}]\, d\xi,
			\end{align*}	
			which gives us the following estimate:
			\begin{align*}
				&\mE	[\|\nab(u(t) - u(s))\|^2_{L^2}] + 2\nu \int_{s}^t \mE[\|\Delta (u(\xi) - u(s))\|^2_{L^2}]\, d\xi\\\nonumber
				&\leq C\Bigl\{\sup_{s\in [0,T]} \mE[\|\Delta u(s)\|^2_{L^2}] + \sup_{\xi\in [0,T]}\mE[\|\nab\cdot(u(\xi)\nab c(\xi))\|^2_{L^2}] \\\nonumber
				&\qquad + \sup_{\xi\in [0,T]}  \mE[\|\nab\cdot(B\nab u(\xi))\|^2_{L^2}] + \sup_{\xi\in [0,T]}\mE[\|\nab(\vb\cdot\nab u(\xi))\|^2_{L^2}]\Bigr\}|t-s|\\\nonumber
				&\leq C\Bigl\{\sup_{s\in [0,T]} \mE[\|\Delta u(s)\|^2_{L^2}] + \sup_{\xi\in [0,T]}\mE[\|u(\xi)\|^4_{H^1}] \\\nonumber
				&\qquad + \sup_{\xi\in [0,T]}  \mE[\|\nab\cdot(B\nab u(\xi))\|^2_{L^2}] + \sup_{\xi\in [0,T]}\mE[\|\nab(\vb\cdot\nab u(\xi))\|^2_{L^2}]\Bigr\}|t-s|.
			\end{align*}
			The proof is complete.
			\item[(c)] Similarly to Part (b), applying the It\^o formula to $\Phi(u(t)) = \|\nab(u(t) - u(s))\|^4_{L^2}$, we obtain
			\begin{align*}
				&\|\nab(u(t) - u(s))\|^4_{L^2} \\\nonumber
				&= 4\nu \int_{s}^t \|\nab(u(\xi) - u(s))\|^2_{L^2}\bigl(\nab(u(\xi) - u(s)),\nab(\Delta u(\xi))\bigr)\, d\xi\\\nonumber
				&\qquad-4\chi\int_{s}^t \|\nab(u(\xi) - u(s))\|^2_{L^2}(\nab(u(\xi) - u(s)),\nab(\nab\cdot(u(\xi)\nab c(\xi))))\, d\xi\\\nonumber
				&\qquad + 4\delta \int_{s}^t  \|\nab(u(\xi) - u(s))\|^2_{L^2}(\nab(u(\xi) - u(s)),\nab(\vb\cdot\nab u(\xi)))\, dW(\xi)\\\nonumber
				&\qquad+2\delta^2 \int_{s}^t  \|\nab(u(\xi) - u(s))\|^2_{L^2}(\nab(u(\xi) - u(s)),\nab(\nab\cdot(B\nab u(\xi))))\, d\xi\\\nonumber
				&\qquad+ 2\delta^2 \int_s^t\|\nab(u(\xi) - u(s))\|^2_{L^2} \|\nab(\vb\cdot\nab u(\xi))\|^2_{L^2}\, d\xi\\\nonumber
				&\qquad+ 4\delta^2\int_{s}^t \bigl(\nab(u(\xi)- u(s)), \nab(\vb\cdot\nab u(\xi))\bigr)^2\, d\xi.
			\end{align*}
			
			Taking the expectation and using integration by parts on the above identity, we obtain
			\begin{align*}
				&\mE[\|\nab(u(t) - u(s))\|^4_{L^2} + 4\nu \int_{s}^t\mE[\|\nab(u(\xi) - u(s))\|^2_{L^2}\|\Delta(u(\xi) - u(s))\|^2_{L^2}]\, d\xi \\\nonumber
				&= -4\nu \int_{s}^t \mE[\|\nab(u(\xi) - u(s))\|^2_{L^2}\bigl(\Delta(u(\xi) - u(s)),\Delta u(s)\bigr)]\, d\xi\\\nonumber
				&\qquad+4\chi\int_{s}^t\mE[\|\nab(u(\xi) - u(s))\|^2_{L^2}(\Delta(u(\xi) - u(s)),\nab\cdot(u(\xi)\nab c(\xi)))]\, d\xi\\\nonumber
				&\qquad-2\delta^2 \int_{s}^t \mE[\|\nab(u(\xi) - u(s))\|^2_{L^2}(\Delta(u(\xi) - u(s)),\nab\cdot(B\nab u(\xi)))]\, d\xi\\\nonumber
				&\qquad+ 2\delta^2 \int_s^t\mE[\|\nab(u(\xi) - u(s))\|^2_{L^2} \|\nab(\vb\cdot\nab u(\xi))\|^2_{L^2}]\, d\xi\\\nonumber
				&\qquad+ 4\delta^2\int_{s}^t \mE[\bigl(\nab(u(\xi)- u(s)), \nab(\vb\cdot\nab u(\xi))\bigr)^2]\, d\xi\\\nonumber
				&:=I_1 + I_2 + I_3 + I_4 + I_5.
			\end{align*}
			
			Next, using the Cauchy-Schwarz inequality and using Lemma \ref{Stability_PDE} and Lemma \ref{Lemma_Holder} part (b) we have
			\begin{align*}
				I_1 &\leq \nu \int_{s}^t\mE[\|\nab(u(\xi) - u(s))\|^2_{L^2}\|\Delta(u(\xi) - u(s))\|^2_{L^2}]\, d\xi\\\nonumber
				&\qquad + C \int_{s}^t\mE[\|\nab(u(\xi) - u(s))\|^2_{L^2}\|\Delta u(s)\|^2_{L^2}]\, d\xi\\\nonumber
				&\leq \nu \int_{s}^t\mE[\|\nab(u(\xi) - u(s))\|^2_{L^2}\|\Delta(u(\xi) - u(s))\|^2_{L^2}]\, d\xi\\\nonumber
				&\qquad + C\int_{s}^t\mE[\|\nab(u(\xi) - u(s))\|^2_{L^2}]\, d\xi\\\nonumber
				&\leq \nu \int_{s}^t\mE[\|\nab(u(\xi) - u(s))\|^2_{L^2}\|\Delta(u(\xi) - u(s))\|^2_{L^2}]\, d\xi + C|t-s|^2.
			\end{align*}
			
			Similarly,
			\begin{align*}
				I_2 &\leq \nu \int_{s}^t\mE[\|\nab(u(\xi) - u(s))\|^2_{L^2}\|\Delta(u(\xi) - u(s))\|^2_{L^2}]\, d\xi\\\nonumber
				&\qquad + C \int_{s}^t\mE[\|\nab(u(\xi) - u(s))\|^2_{L^2}\|\nab\cdot(u(\xi)\nab c(\xi))\|^2_{L^2}]\, d\xi\\\nonumber
				&\leq \nu \int_{s}^t\mE[\|\nab(u(\xi) - u(s))\|^2_{L^2}\|\Delta(u(\xi) - u(s))\|^2_{L^2}]\, d\xi\\\nonumber
				&\qquad + C \int_{s}^t\mE[\|\nab(u(\xi) - u(s))\|^2_{L^2}]\, d\xi\\\nonumber
				&\leq  \nu \int_{s}^t\mE[\|\nab(u(\xi) - u(s))\|^2_{L^2}\|\Delta(u(\xi) - u(s))\|^2_{L^2}]\, d\xi + C|t-s|^2.\\\nonumber
				I_3 &\leq  \nu \int_{s}^t\mE[\|\nab(u(\xi) - u(s))\|^2_{L^2}\|\Delta(u(\xi) - u(s))\|^2_{L^2}]\, d\xi + C|t-s|^2.
			\end{align*}
			
			Next using Lemma \ref{Lemma_Holder} and Lemma \ref{Stability_PDE} (d) with the deterministic intial data $u_0$, we obtain
			\begin{align*}
				I_4 &\leq C\int_{s}^t \mE[\|\nab(u(\xi) - u(s))\|^2_{L^2}]\, d\xi \leq CC_2|t-s|^2.\\\nonumber
				I_5 &\leq C\int_{s}^t \mE[\|\nab(u(\xi) - u(s))\|^2_{L^2}]\, d\xi \leq CC_2|t-s|^2.
			\end{align*}
			
			Collecting all the estimates from $I_1, \cdots, I_5$ and substituting them into the original identity above, we obtain the desired estimate. The proof is complete.

		\end{enumerate}
	\end{proof}
	
	\section{Semi discretization in time}\label{sec3}
	\subsection{Semi-implicit Crank-Nicolson method}
	In this section, we present the time discretization of \eqref{Weak_formulation} by using the Crank-Nicolson method, which is a natural method to approximate the Stratonovich noise. 
	
	Let $k = \frac{T}{M}$ be the uniform time step of a partition $\{t_m\}^M_{m=0}$ over the time interval $[0,T]$. Let's define the following notations for our numerical scheme later.
	\begin{align*}
		\phi^{m+\frac12} = \frac12\bigl(\phi^{m+1} + \phi^{m}\bigr), \qquad \phi^{m-\frac12} = \frac12\bigl(\phi^{m} + \phi^{m-1}\bigr), \qquad \phi^{\frac12} = \frac12\bigl(\phi^1 + \phi^0\bigr).
	\end{align*}

	Then we can state the Crank-Nicolson method for \eqref{Weak_formulation} in the following algorithm.
	
	\medskip
	{\bf Algorithm 1:} Let $u^0 =  u_0$. For any $0\leq m \leq M-1$, we do following iterations:
	Compute $\bigl(u^{m+1}, \vsigma^{m+1}, c^{m+1}\bigr) \in H^1_{per}(D)\times \vH^1_{per}(D)\times L^2_{per}(D)$ such that $\mP$- a.s.
	\begin{align}
		\label{eq_semi_discrete_sigma}	&\bigl(\vsigma^{m+1},\vvarphi\bigr) + \bigl(\nab\cdot\vsigma^{m+1},\nab\cdot\vvarphi\bigr) + \bigl(rot\, \vsigma^{m+1},rot\, \vvarphi\bigr) = -\bigl(u^{m}, \nab\cdot\vvarphi\bigr),\\
		\medskip
		\label{eq_semi_discrete_u}	&\bigl(u^{m+1}- u^{m},v\bigr) + k \nu\bigl(\nab u^{m+\frac12},\nab v\bigr) =  \chi k\bigl(u^{m+\frac12}\vsigma^{m+\frac12},\nab v\bigr) \\\nonumber
		&\qquad\qquad\qquad\qquad\qquad\qquad\qquad\qquad\qquad+ \delta\bigl(\vb\cdot\nab u^{m+\frac12}\Delta W_m, v\bigr),\\\medskip
		\label{eq_semi_discrete_c}	&\bigl(c^{m+1}, \psi\bigr) = \bigl(\nab\cdot \vsigma^{n+1}, \psi\bigr)+ \bigl(u^{m+1},\psi\bigr),
	\end{align}
	for all $\bigl(v, \vvarphi,\psi\bigr) \in H^1_{per}(D)\times \vH^1_{per}(D)\times L^2_{per}(D)$, and $\Delta W_m = W(t_{m+1}) - W(t_m)$.
	
	\medskip
	Next, we present the stability estimates of Algorithm 1. 
	\begin{lemma}\label{Lemma_stability} Let $(u^m, \vsigma^m,c^m)$ be the solution of Algorithm 1. Let $u_0$ such that  
		\begin{align}\label{stability_criteria}
			1 - \frac{32\chi^2C_L^4T\kappa_1}{\nu} >0,
		\end{align}
		
		where $\kappa_1 = \|u_0\|^2_{L^2} + \frac{3\chi^2 C_L^4 k}{\nu}\|u_0\|^4_{L^2}$.
		
		Then, they hold that
		\begin{enumerate}
			\item[(a)] If $u_0\in L_{per}^2(D)$, then there exists $K_1 = K_1(u_0)$ such that
			\begin{align}\label{eq3.4}
				\max_{1\leq m \leq M}\|u^m\|^2_{L^2} + \nu k\sum_{m=1}^{M-1}\|\nab u^{m+\frac12}\|^2_{L^2} &\leq K_1,\qquad\mP-a.s.,\\
				\label{eq3.5}	\max_{1\leq m \leq M}\|\vsigma^m\|^2_{H^1} &\leq K_1,\qquad\mP-a.s.,\\
				\label{eq3.6}			\max_{1\leq m \leq M} \|c^m\|^2_{L^2}  &\leq K_1\qquad\mP-a.s..
			\end{align}
			\item[(b)] If $u_0\in H_{per}^1(D)$, then there exists $K_2 = K_2(u_0,K_1)$ such that
			\begin{align}
				\label{eq3.7}	\max_{1\leq m \leq M}\|\nab u^m\|^2_{L^2} + \nu k\sum_{m=1}^{M-1}\|\Delta u^{m+\frac12}\|^2_{L^2} &\leq K_2,\qquad\mP-a.s.,\\
				\label{eq3.8}	\max_{1\leq m \leq M}\left(\|\nab \vsigma^m\|^2_{L^2} + \|\nab(\nab\cdot\vsigma)^m\|^2_{L^2}\right) &\leq K_2,\qquad\mP-a.s..
			\end{align}
			\item[(c)] If $u_0\in H_{per}^2(D)$, then there exists $K_3 = K_3(u_0,K_2)>0$ such that
			\begin{align}
				\label{eq3.9}	\max_{1\leq m \leq M}\|\Delta u^m\|^2_{L^2} + \nu k\sum_{m=1}^{M-1}\|\nab \Delta u^{m+\frac12}\|^2_{L^2} &\leq K_3,\qquad\mP-a.s..
			\end{align}
		\end{enumerate}
		
	\end{lemma}
	\begin{proof}
		\begin{enumerate}
			\item[(a)] 	First, to prove \eqref{eq3.4}, taking $v= u^{m+\frac12}$ in \eqref{eq_semi_discrete_u} and use the identity $(a-b)(a+b) = a^2 - b^2$ we obtain
			\begin{align}\label{eq_3.6}
				&\frac{1}{2}\bigl[\|u^{m+1}\|^2_{L^2} - \|u^{m}\|^2_{L^2} \bigr] + \nu k\|\nab u^{m+\frac12}\|^2_{L^2}\\\nonumber
				& = \chi k\bigl(u^{m+\frac12}\vsigma^{m+\frac12}, \nab u^{m+1}\bigr) + \delta\bigl(\vb\cdot\nab u^{m+\frac12}\Delta W_m,  u^{m+\frac12}\bigr)\\\nonumber
				&:= I^a_1 + I^a_2.
			\end{align}
			
			Using integration by parts, we get $I^a_2 =0$. Next, using integration by parts and Cauchy-Schwarz's inequality and then Ladyzhenskaya's inequality, we obtain
			\begin{align}\label{eq_3.7}
				I^a_1&=-\chi k \bigl((u^{m+\frac12})^2, \nab\cdot\vsigma^{m+\frac12}\bigr)\\\nonumber
				&\leq \chi k \|u^{m+\frac12}\|^2_{L^4}\|\nab\cdot \vsigma^{m+\frac12}\|_{L^2}\\\nonumber
				&\leq \chi k C_L^2 \|u^{m+\frac12}\|_{L^2}\|\nab u^{m+\frac12}\|_{L^2}\|\nab\cdot \vsigma^{m+\frac12}\|_{L^2}
			\end{align}
			
			Next, choosing $\vvarphi = \vsigma^{m+1}$ in \eqref{eq_semi_discrete_sigma}, we obtain
			\begin{align*}
				\|\vsigma^{m+1}\|^2_{H^1} \leq \|u^{m}\|^2_{L^2},
			\end{align*}
			which together with \eqref{eq_3.7} imply that 
			\begin{align}\label{eq_3.8}
				I^a_1 &\leq  \chi k C_L^2 \|u^{m+\frac12}\|_{L^2}\|\nab u^{m+\frac12}\|_{L^2}\|u^{m-\frac12}\|_{L^2}\\\nonumber
				&\leq \frac{\nu k}{4} \|\nab u^{m+\frac12}\|^2_{L^2} + \frac{\chi^2C_L^4k}{\nu}\|u^{m+\frac12}\|^2_{L^2}\|u^{m-\frac12}\|^2_{L^2}\\\nonumber
				&\leq \frac{\nu k}{4} \|\nab u^{m+\frac12}\|^2_{L^2} + \frac{\chi^2C_L^4k}{2\nu}\bigl(\|u^{m+\frac12}\|^4_{L^2} + \|u^{m-\frac12}\|^4_{L^2}\bigr).
			\end{align}
			
			Substituting \eqref{eq_3.8} into \eqref{eq_3.6}, we arrive at
			\begin{align}\label{eq_3.9}
				\|u^{m+1}\|^2_{L^2} - \|u^{m}\|^2_{L^2} + \frac{3\nu k}{2}\|\nab u^{m+\frac12}\|^2_{L^2}
				&\leq \frac{\chi^2C_L^4k}{\nu}\bigl(\|u^{m+\frac12}\|^4_{L^2} + \|u^{m-\frac12}\|^4_{L^2}\bigr).
			\end{align}
			
			Next, applying the summation $\sum_{m=1}^{\ell}$ for any $1\leq \ell \leq M-1$, we obtain
			\begin{align}\label{eq_3.13}
				&\|u^{\ell+1}\|^2_{L^2} + \frac{3\nu k}{2}\sum_{m=1}^{\ell} \|\nab u^{m+\frac12}\|^2_{L^2} \\\nonumber
				&\leq \|u^1\|^2_{L^2} + \frac{\chi^2C_L^4k}{\nu}\sum_{m=1}^{\ell}\bigl(\|u^{m+\frac12}\|^4_{L^2} + \|u^{m-\frac12}\|^4_{L^2}\bigr)\\\nonumber
			\end{align}
			
			In addition, using $(a+b)^4 \leq 2^3(a^4 + b^4)$ to the last term on the right-hand side of \eqref{eq_3.13}, we obtain
			\begin{align}\label{eq_3.14}
				&\|u^{\ell+1}\|^2_{L^2} +\nu k\sum_{m=1}^{\ell} \|\nab u^{m+\frac12}\|^2_{L^2} \\\nonumber
				&\leq \|u^1\|^2_{L^2} + \frac{8\chi^2C_L^4k}{\nu}\sum_{m=1}^{\ell}\bigl(\|u^{m+1}\|^4_{L^2} + 2\|u^m\|^4_{L^2} + \|u^{m-1}\|^4_{L^2}\bigr)\\\nonumber
				&\leq \|u^1\|^2_{L^2} + \frac{32\chi^2C_L^4k}{\nu}\sum_{m=1}^{\ell}\|u^{m+1}\|^4_{L^2}.
			\end{align}
			
			Next, taking $v= u^{\frac12}$ in \eqref{eq_semi_discrete_u} with $m=1$, we obtain
			\begin{align*}
				\frac12\bigl[\|u^1\|^2_{L^2} - \|u^0\|^2_{L^2}\bigr] + \nu k\|\nab u^{\frac12}\|^2_{L^2} &= \chi k \bigl(u^{\frac12}\vsigma^{\frac12}, \nab u^{\frac12}\bigr)\\\nonumber
				&=-\chi k \bigl((u^{\frac12})^2,\nab\cdot\vsigma^{\frac12}\bigr)\\\nonumber
				&\leq \chi k C_L^2 \|u^{\frac12}\|_{L^2}\|\nab u^{\frac12}\|_{L^2}\|\nab\cdot\vsigma^{\frac12}\|_{L^2}
			\end{align*}
			which implies that
			\begin{align*}
				\frac12\|u^1\|^2_{L^2} + \frac{3\nu k}{4}\|\nab u^{\frac12}\|^2_{L^2}& \leq \frac12\|u^0\|^2_{L^2} + \frac{\chi^2 C_L^4 k}{\nu} \|u^{\frac12}\|^2_{L^2}\|\nab\cdot \vsigma^{\frac12}\|^2_{L^2}\\\nonumber
				&\leq \frac12\|u^0\|^2_{L^2} + \frac{\chi^2 C_L^4 k}{2\nu}\bigl( \|u^{\frac12}\|^4_{L^2} +\|\nab\cdot \vsigma^{\frac12}\|^4_{L^2}\bigr).
			\end{align*}
			
			So, 
			\begin{align}\label{eq_3.16}
				\|u^1\|^2_{L^2} + \frac{3\nu k}{2}\|\nab u^{\frac12}\|^2_{L^2}& \leq \|u^0\|^2_{L^2} + \frac{\chi^2 C_L^4 k}{\nu} \|u^{1}\|^4_{L^2} + \frac{\chi^2 C_L^4 k}{\nu} \|u^{0}\|^4_{L^2}\\\nonumber
				&\qquad + \frac{\chi^2 C_L^4 k}{\nu} \|\nab\cdot\vsigma^{1}\|^4_{L^2} + \frac{\chi^2 C_L^4 k}{\nu} \|\nab\cdot\vsigma^{0}\|^4_{L^2}\\\nonumber
				& \leq \|u^0\|^2_{L^2} + \frac{\chi^2 C_L^4 k}{\nu} \|u^{1}\|^4_{L^2} + \frac{2\chi^2 C_L^4 k}{\nu} \|u^{0}\|^4_{L^2}\\\nonumber
				&\qquad  + \frac{\chi^2 C_L^4 k}{\nu} \|\nab\cdot\vsigma^{0}\|^4_{L^2} \\\nonumber
				& \leq \|u^0\|^2_{L^2} + \frac{\chi^2 C_L^4 k}{\nu} \|u^{1}\|^4_{L^2} + \frac{2\chi^2 C_L^4 k}{\nu} \|u^{0}\|^4_{L^2}\\\nonumber
				&\qquad  + \frac{\chi^2 C_L^4 k}{\nu} \|c_0\|^4_{H^2} \\\nonumber
				&\leq \|u^0\|^2_{L^2} + \frac{\chi^2 C_L^4 k}{\nu} \|u^{1}\|^4_{L^2} + \frac{2\chi^2 C_L^4 k}{\nu} \|u^{0}\|^4_{L^2}\\\nonumber
				&\qquad  + \frac{\chi^2 C_L^4 k}{\nu} \|u_0\|^4_{L^2} \\\nonumber
				&:= \kappa_1(u_0) + \frac{\chi^2 C_L^4 k}{\nu} \|u^{1}\|^4_{L^2} 
			\end{align}
			
			Next, substituting \eqref{eq_3.16} into \eqref{eq_3.14} we arrive at
			\begin{align}\label{eq_3.18}
				\|u^{\ell+1}\|^2_{L^2} +\nu k\sum_{m=1}^{\ell} \|\nab u^{m+\frac12}\|^2_{L^2} 
				&\leq \kappa_1(u_0)+ \frac{32\chi^2C_L^4}{\nu} k\sum_{m=0}^{\ell} \|u^{m+1}\|^4_{L^2}.
			\end{align}
			
			Using Lemma \ref{nonlinear_Gronwall} with $a =\kappa_1(u_0), b(n) =0, k(n) = \frac{32\chi^2C_L^4 k}{\nu}$, and $p=2$ we obtain
			\begin{align}\label{eq_3.10}
				\|u^{\ell+1}\|^2_{L^2} + \nu k\sum_{m=1}^{\ell} \|\nab u^{m+\frac12}\|^2_{L^2} 
				&\leq \frac{1}{{\frac{1}{\kappa_1} - \frac{32 \chi^2C_L^4t_{\ell}}{\nu}}} \leq \frac{\kappa_1}{{1 - \frac{32\chi^2C_L^4T\kappa_1}{\nu}}}:=K_1,
			\end{align}
			
			where $1 - \frac{32\chi^2C_L^4T\kappa}{\nu} >0$ by using the assumption \eqref{stability_criteria} of $u_0$.
			
			Finally, the desired estimates are established by taking the maximum over $0\leq \ell\leq M-1$ to \eqref{eq_3.10}.
			
			Lastly, it is clear that $\eqref{eq3.5}$ is obtained by taking $\vvarphi= \vsigma^{m+1}$ and using $\eqref{eq3.4}$. Similarly, $\eqref{eq3.6}$ is a consequence of $\eqref{eq3.4}$, and $\eqref{eq3.5}$.\\

			\item[(b)] 	In the next steps, we are going to prove $\eqref{eq3.7}$. Taking $v = -\Delta u^{m+\frac12}$ in \eqref{eq_semi_discrete_u} and using integration by parts, we obtain
			\begin{align}\label{Delta_stability}
				&\frac12\bigl[\|\nab u^{m+1}\|^2_{L^2} - \|\nab u^{m}\|^2_{L^2}\bigr] + \nu k \|\Delta u^{m+\frac12}\|^2_{L^2}\\\nonumber
				&=\chi k \bigl(\nab[u^{m+\frac12}\vsigma^{m+\frac12}],\Delta u^{m+\frac12}\bigr) - \delta \bigl(\vb\cdot\nab u^{m+\frac12},\Delta u^{m+\frac12}\bigr)\Delta W_m\\\nonumber
				&:= I^b_1 + I^b_2,
			\end{align}
			where $I^b_2 =0$ by using integration by parts. In addition, using Cauchy-Schwarz's inequality and then Ladyzhenskaya's inequality, we have
			\begin{align*}
				I^b_1 &= \chi k\bigl(\nab u^{m+\frac12}\vsigma^{m+\frac12} + u^{m+\frac12}\nab\cdot \vsigma^{m+\frac12},\Delta u^{m+\frac12}\bigr)\\\nonumber
				&\leq \frac{2\chi^2k}{\nu}\Bigl[\|\nab u^{m+\frac12}\vsigma^{m+\frac12}\|^2_{L^2} + \| u^{m+\frac12}\nab\cdot\vsigma^{m+\frac12}\|^2_{L^2}\Bigr] + \frac{\nu k}{4} \|\Delta u^{m+\frac12}\|^2_{L^2}\\\nonumber
				&\leq  \frac{2\chi^2k}{\nu}\Bigl[\|\nab u^{m+\frac12}\|^2_{L^4}\|\vsigma^{m+\frac12}\|^2_{L^4} + \| u^{m+\frac12}\|^2_{L^4}\|\nab\cdot\vsigma^{m+\frac12}\|^2_{L^4}\Bigr] + \frac{\nu k}{4} \|\Delta u^{m+\frac12}\|^2_{L^2}\\\nonumber
				&\leq  \frac{2\chi^2k}{\nu}C_L^2\Bigl[\|\nab u^{m+\frac12}\|_{L^2}\|\Delta u^{m+\frac12}\|_{L^2}\|\vsigma^{m+\frac12}\|_{L^2}\|\vsigma^{m+\frac12}\|_{H^1} \\\nonumber
				&\qquad+ \| u^{m+\frac12}\|_{L^2}\|\nab u^{m+\frac12}\|_{L^2}\|\nab\cdot\vsigma^{m+\frac12}\|_{L^2}\|\nab(\nab\cdot\vsigma^{m+\frac12})\|_{L^2}\Bigr] + \frac{\nu k}{4} \|\Delta u^{m+\frac12}\|^2_{L^2}\\\nonumber
				&\leq  \frac{2\chi^2k}{\nu}C_L^2K_1^2\Bigl[\|\nab u^{m+\frac12}\|_{L^2}\|\Delta u^{m+\frac12}\|_{L^2} \\\nonumber&\qquad+ \|\nab u^{m+\frac12}\|_{L^2}\|\nab (\nab\cdot\vsigma^{m+\frac12})\|_{L^2}\Bigr] + \frac{\nu k}{4} \|\Delta u^{m+\frac12}\|^2_{L^2}\\\nonumber
				&\leq \frac{\nu k}{4} \|\Delta u^{m+\frac12}\|^2_{L^2} + \frac{\nu k}{4} \|\Delta u^{m+\frac12}\|^2_{L^2} + \frac{4\chi^4}{\nu^3} C_L^4 K_1^4k\|\nab u^{m+\frac12}\|^2_{L^2} \\\nonumber
				&\qquad +\frac{2\chi^2 k}{\nu}C_L^2 K_1^2 \|\nab u^{m+\frac12}\|_{L^2}\|\nab(\nab\cdot\vsigma^{m+\frac12})\|_{L^2}.
			\end{align*}
			
			Next, taking $\vvarphi = -\nab (\nab\cdot \sigma^{m+\frac12})$ in \eqref{eq_semi_discrete_sigma} and using integration by parts, we obtain
			\begin{align}\label{eq_3.20}
				\|\nab\cdot \vsigma^{m+\frac12}\|^2_{L^2} + \|\nab(\nab\cdot\vsigma^{m+\frac12})\|^2_{L^2} \leq \|\nab u^{m-\frac12}\|^2_{L^2}.
			\end{align}
			Replacing \eqref{eq_3.20} to the right-hand side of $I^b_1$, we obtain
			\begin{align}\label{eq_3.21}
				I^b_1 &\leq \frac{\nu k}{4} \|\Delta u^{m+\frac12}\|^2_{L^2} + \frac{\nu k}{4} \|\Delta u^{m+\frac12}\|^2_{L^2} + \frac{4\chi^4}{\nu^3} C_L^4 K_1^4k\|\nab u^{m+\frac12}\|^2_{L^2} \\\nonumber
				&\qquad +\frac{2\chi^2 k}{\nu}C_L^2 K_1^2 \|\nab u^{m+\frac12}\|_{L^2} \|\nab u^{m-\frac12}\|_{L^2}\\\nonumber
				&\leq \frac{\nu k}{4} \|\Delta u^{m+\frac12}\|^2_{L^2} + \frac{\nu k}{4} \|\Delta u^{m+\frac12}\|^2_{L^2} + \frac{4\chi^4}{\nu^3} C_L^4 K_1^4k\|\nab u^{m+\frac12}\|^2_{L^2} \\\nonumber
				&\qquad +\frac{\chi^2 k}{\nu}C_L^2 K_1^2 \|\nab u^{m+\frac12}\|^2_{L^2} + \frac{\chi^2 k}{\nu}C_L^2 K_1^2\|\nab u^{m-\frac12}\|^2_{L^2}\\\nonumber
				&= \frac{\nu k}{2} \|\Delta u^{m+\frac12}\|^2_{L^2} 
				+ \Bigl(1+ \frac{4\chi^2}{\nu^2} C_L^2 K_1^2\Bigr)\frac{\chi^2}{\nu} C_L^2 K_1^2k\|\nab u^{m+\frac12}\|^2_{L^2}\\\nonumber
				&\qquad+ \frac{\chi^2 k}{\nu}C_L^2 K_1^2\|\nab u^{m-\frac12}\|^2_{L^2}.
			\end{align}
			Substituting \eqref{eq_3.21} into \eqref{Delta_stability} and then applying $\sum_{m=1}^{\ell}$, for $1\leq \ell\leq M-1$, we obtain
			\begin{align}\label{eq_3.233}
				&	\frac12\|\nab u^{m+1}\|^2_{L^2} +\frac{\nu k}{2}\sum_{m=1}^{\ell}\|\Delta u^{m+\frac12}\|^2_{L^2} \\\nonumber
				&\leq \Bigl(1+ \frac{2\chi^2}{\nu^2} C_L^2 K_1^2\Bigr)\frac{2\chi^2}{\nu} C_L^2 K_1^2k\sum_{m=1}^{\ell}\|\nab u^{m+\frac12}\|^2_{L^2} + \frac12\|\nab u^1\|^2_{L^2}.
			\end{align}
			
			To estimate $\|\nab u^1\|^2_{L^2}$, taking $m=0$ in \eqref{Delta_stability}, we obtain
			\begin{align}\label{eq3.23}
				\frac12	\|\nab u^1\|^2_{L^2} + \nu k \|\Delta u^{\frac12}\|^2_{L^2}  &= \frac12\|\nab u_0\|^2_{L^2}+ \chi k\bigl(\nab[u^{\frac12}\vsigma^{\frac12}], \Delta u^{\frac12}\bigr).
			\end{align}
			
			To bound the last term of \eqref{eq3.23}, we proceed as $I^b_1$ with $m=0$ to obtain
			\begin{align}\label{eq3.24}
				\chi k\bigl(\nab[u^{\frac12}\vsigma^{\frac12}], \Delta u^{\frac12}\bigr) &\leq \frac{\nu k}{2} \|\Delta u^{\frac12}\|^2_{L^2} + \frac{4\chi^4}{\nu^3} C_L^4 K_1^4k\|\nab u^{\frac12}\|^2_{L^2} \\\nonumber
				&\qquad +\frac{2\chi^2 k}{\nu}C_L^2 K_1^2 \|\nab u^{\frac12}\|_{L^2}\|\nab(\nab\cdot\vsigma^{\frac12})\|_{L^2}.
			\end{align}
			
			Next, we notice that using \eqref{eq_semi_discrete_sigma} with $m=0$ and $\varphi = -\nab(\nab\cdot \vsigma^1)$, we can control 
			\begin{align*}
				\|\nab(\nab\cdot\vsigma^{\frac12})\|_{L^2} & =\frac12\|\nab(\nab\cdot \vsigma^1) +\nab(\nab\cdot \vsigma^0)  \|_{L^2} \\\nonumber
				&\leq \frac12\|\nab(\nab\cdot \vsigma^1) \|_{L^2} + \frac12\|\nab(\nab\cdot \vsigma^0) \|_{L^2}\\\nonumber
				& \leq \frac{1}{2}\|\nab u^0\|_{L^2} + \frac12\|u_0\|_{H^1}.
			\end{align*}
			
			With this, \eqref{eq3.24} and \eqref{eq3.23}, we obtain
			\begin{align*}
				\frac12	\|\nab u^1\|^2_{L^2} + \frac{\nu k}{2} \|\Delta u^{\frac12}\|^2_{L^2}  &\leq \frac12\|\nab u_0\|^2_{L^2}+  \frac{4\chi^4}{\nu^3} C_L^4 K_1^4k\|\nab u^{\frac12}\|^2_{L^2} \\\nonumber
				&\qquad +\frac{\chi^2 k}{\nu}C_L^2 K_1^2 \|\nab u^{\frac12}\|^2_{L^2} +\frac{\chi^2 k}{\nu}C_L^2 K_1^2\| u^0\|^2_{H^1} .
			\end{align*}
			
			Substituting this into \eqref{eq_3.233}, we arrive at
			\begin{align}
				&	\frac12\|\nab u^{m+1}\|^2_{L^2} +\frac{\nu k}{2}\sum_{m=1}^{\ell}\|\Delta u^{m+\frac12}\|^2_{L^2} \\\nonumber
				&\leq \Bigl(1+ \frac{2\chi^2}{\nu^2} C_L^2 K_1^2\Bigr)\frac{2\chi^2}{\nu} C_L^2 K_1^2k\sum_{m=0}^{\ell}\|\nab u^{m+\frac12}\|^2_{L^2} + \frac12\|\nab u_0\|^2_{L^2}\\\nonumber
				&\qquad+\frac{\chi^2 k}{\nu}C_L^2 K_1^2\|u^0\|^2_{H^1}.
			\end{align}

			Using the standard discrete Gronwall inequality, we obtain
			\begin{align}\label{eq_3.23}
				&\|\nab u^{\ell+1}\|^2_{L^2} +\frac{\nu k}{2}\sum_{m=1}^{\ell}\|\Delta u^{m+\frac12}\|^2_{L^2} \\\nonumber
				&\qquad\leq \exp\Bigl[\Bigl(1+ \frac{4\chi^2}{\nu^2} C_L^2 K_1^2\Bigr)\frac{2\chi^2}{\nu} C_L^2 K_1^2T\Bigr]\kappa_2:=K_2,
			\end{align}
			where $\kappa_2 = \|\nab u_0\|^2_{L^2} + \frac{2\chi^2 k}{\nu}C_L^2 K_1^2\|u^0\|^2_{H^1}.$
			
			Finally, the desired estimates are established by taking the maximum over $0\leq \ell\leq M-1$ to \eqref{eq_3.23}. Lastly, it is clear that $(b)_2$ is obtained by taking $\vvarphi =\nab(\nab\cdot \sigma^m)$ in \eqref{eq_semi_discrete_sigma} and using $(b)_1$. 
			\item[(c)] Similarly to part (b), we prove \eqref{eq3.9} by taking $v = \Delta^2 u^{m+\frac12}$ in \eqref{eq_semi_discrete_u} and using integration by parts twice, we get
			\begin{align}\label{eq3.27}
				&	\frac12[\|\Delta u^{m+1}\|^2_{L^2} - \|\Delta u^m\|^2_{L^2}] + \nu k \|\nab(\Delta u^{m+\frac12})\|^2_{L^2} &\\\nonumber
				&= \chi k \bigl(\nab(\nab\cdot(u^{m+\frac12}\vsigma^{m+\frac12})), \nab(\Delta u^{m+\frac12})\bigr) + 0: = I_1^c.
			\end{align}
			
			To estimate $I_1^c$ of \eqref{eq3.27}, we follow the same lines as the estimation of $I^b_1$ in part (b) and then use the standard Gronwall inequality to complete the proof. So, we skip the detailed proof to save space.
			
		\end{enumerate}

	\end{proof}

	\subsection{Semi-discrete error estimates}
	Next, we state and prove error estimates of Algorithm 1.
	\begin{theorem}\label{Thm_error_semi} Let $(u^m,\vsigma^m,c^m)$ be the solution from Algorithm 1. Suppose that $u_0 \in L^{\infty}(\Omega; H^2_{per}(D))$ and assume that the condition \eqref{stability_criteria} of Lemma \ref{Lemma_stability} satisfies. Additionally, assume that for given $\nu>0, \chi>0$ and $\vb\in \mathbb{R}^2$ 
		\begin{align}\label{assume_delta}
			\delta^2 < \min\left\{\frac{7\nu}{8C(D)|\vb|^2}, \frac{15\nu}{16|\vb|^2}\right\}.
		\end{align} 
		
		Then, there exists a constant $C = C(C_1,C_2,K_1, K_2,K_3)>0$ such that
		\begin{align}
			\label{estimate_u}	\max_{1\leq m \leq M}\mE\bigl[\|u(t_m) - u^m\|^2_{L^2}\bigr] + \mE\Bigl[\nu k\sum_{m=1}^M \|\nab (u(t_m) - u^m)\|^2_{L^2}\Bigr] &\leq Ck,\\
			\label{estimate_sigma}	\max_{1\leq m \leq M}\mE\bigl[\|\vsigma(t_m) - \vsigma^m\|^2_{H^1}\bigr] &\leq Ck,\\
			\label{estimate_c}		\max_{1\leq m \leq M}\mE\bigl[\|c(t_m) - c^m\|^2_{L^2}\bigr] &\leq Ck.
		\end{align}
	\end{theorem}
	\begin{proof} First, \eqref{estimate_sigma} and \eqref{estimate_c} can be obtained directly from \eqref{estimate_u}. So, we just need to give the proof of \eqref{estimate_u}.
		
		Denote $e^m_u := u(t_m) - u^m$, $e_{\vsigma}^m := \sigma(t_m) - \vsigma^m$ and $e_c^m := c(t_m) - c^m$. We also denote $u(t_{m+\frac12}) := \frac12\bigl[u(t_{m+1}) + u(t_m)\bigr]$ and $u(t_{m-\frac12}) := \frac12\bigl[u(t_m) + u(t_{m-1})\bigr]$.
		
		Subtract \eqref{eq_u} from \eqref{eq_semi_discrete_u}, we obtain the following error equation.
		\begin{align}\label{eq_3.27}
			&\bigl(e_u^{m+1} - e^m_u,v\bigr) + \nu k\bigl(\nab e_u^{m+\frac12}, \nab v\bigr) \\\nonumber
			&= \nu \int_{t_m}^{t_{m+1}} \bigl(\nab u(t_{m+\frac12}) - \nab u(s), \nab v\bigr)\, ds \\\nonumber
			&+ \chi\int_{t_m}^{t_{m+1}} \bigl(u(s)\vsigma(s) - u^{m+\frac12}\vsigma^{m+\frac12}, \nab v\bigr)\, ds \\\nonumber
			&+ \delta \Bigl(\int_{t_m}^{t_{m+1}}\vb\cdot\nab u(s)\circ dW(s) - \vb\cdot \nab u^{m+\frac12}\Delta W_m, v\Bigr).
		\end{align}
		
		Taking $v = e_u^{m+1}$ in \eqref{eq_3.27}, we obtain the following.
		\begin{align}\label{eq_3.28}
			&\frac12\bigl[\|e^{m+1}_u\|^2_{L^2} - \|e_u^m\|^2_{L^2}\bigr] +\frac12\|e_u^{m+1} - e_u^m\|^2_{L^2}+ \nu k\bigl(\nab e^{m+\frac12},\nab e^{m+1}_u\bigr)\\\nonumber
			&= \nu \int_{t_m}^{t_{m+1}} \bigl(\nab u(t_{m+\frac12}) - \nab u(s), \nab e_u^{m+1}\bigr)\, ds \\\nonumber
			&\qquad+ \chi\int_{t_m}^{t_{m+1}} \bigl(u(s)\vsigma(s) - u^{m+\frac12}\vsigma^{m+\frac12}, \nab e_u^{m+1}\bigr)\, ds \\\nonumber
			&\qquad+ \delta \Bigl(\int_{t_m}^{t_{m+1}}\vb\cdot\nab u(s)\circ dW(s) - \vb\cdot \nab u^{m+\frac12}\Delta W_m, e_u^{m+1}\Bigr).
		\end{align}
		We can analyze the third term on the left side of \eqref{eq_3.28} as follows:
		\begin{align*}
			\nu k\bigl(\nab e^{m+\frac12},\nab e^{m+1}_u\bigr) &=\frac{\nu k}{2}\|\nab e_u^{m+1}\|^2_{L^2} + \frac{\nu k}{2}\|\nab e_u^m\|^2_{L^2} + \frac{\nu k}{2}\bigl(\nab e_u^m,\nab(e_u^{m+1} - e_u^m)\bigr)\\\nonumber
			&=\frac{\nu k}{2}\|\nab e_u^{m+1}\|^2_{L^2} + \frac{\nu k}{2}\|\nab e_u^m\|^2_{L^2} - \frac{\nu k}{2}\bigl(\Delta e_u^m,e_u^{m+1} - e_u^m\bigr).
		\end{align*}
		Substituting this into \eqref{eq_3.28} and rearrange terms, we obtain
		\begin{align}\label{eq_3.29}
			&\frac12\bigl[\|e^{m+1}_u\|^2_{L^2} - \|e_u^m\|^2_{L^2}\bigr] +\frac12\|e_u^{m+1} - e_u^m\|^2_{L^2}  \\\nonumber
			&\qquad\qquad+\frac{\nu k}{2}\|\nab e_u^{m+1}\|^2_{L^2} + \frac{\nu k}{2}\|\nab e_u^m\|^2_{L^2} \\\nonumber
			&= \frac{\nu k}{2}\bigl(\Delta e_u^m,e_u^{m+1} - e_u^m\bigr)
			+\nu \int_{t_m}^{t_{m+1}} \bigl(\nab u(t_{m+\frac12}) - \nab u(s), \nab e_u^{m+1}\bigr)\, ds \\\nonumber
			&\qquad+ \chi\int_{t_m}^{t_{m+1}} \bigl(u(s)\vsigma(s) - u^{m+\frac12}\vsigma^{m+\frac12}, \nab e_u^{m+1}\bigr)\, ds \\\nonumber
			&\qquad+ \delta \Bigl(\int_{t_m}^{t_{m+1}}\vb\cdot\nab u(s)\circ dW(s) - \vb\cdot \nab u^{m+\frac12}\Delta W_m, e_u^{m+1}\Bigr)\\\nonumber
			&:={I + II + III + IV}.
		\end{align}
		
		Using Cauchy-Schwarz's inequality and Lemma \ref{Lemma_stability}, and Lemma \ref{Stability_PDE}, we obtain
		\begin{align*}
			I &\leq \nu^2 k^2 16\|\Delta e_u^{m}\|^2_{L^2} + \frac{1}{64}\|e_u^{m+1} - e_u^m\|^2_{L^2}\\\nonumber
			&\leq C(C_2+K_3) k^2 + \frac{1}{64}\|e_u^{m+1} - e_u^m\|^2_{L^2}.
		\end{align*}
		Using Cauchy-Schwarz's inequality and Lemma \ref{Lemma_Holder}, we obtain
		\begin{align*}
			\mE[II] &\leq 16\nu\int_{t_m}^{t_{m+1}}\mE\bigl[\|\nab(u(t_{m+\frac12}) - u(s))\|^2_{L^2}\bigr]\, ds + \frac{\nu k}{64}\mE\bigl[\|\nab e_u^{m+1}\|^2_{L^2}\bigr]\\\nonumber
			&\leq Ck^{2} + \frac{\nu k}{64}\mE\bigl[\|\nab e_u^{m+1}\|^2_{L^2}\bigr].
		\end{align*}
		
		To control $III$, we add and subtract the term $u(t_{m+\frac12})\vsigma(t_{m+\frac12})$:
		\begin{align}
			III &= \chi \int_{t_m}^{t_{m+1}} \bigl(u(s)\vsigma(s) - u(t_{m+\frac12})\vsigma(t_{m+\frac12}),\nab e_u^{m+1}\bigr)\, ds \\\nonumber
			&\qquad+ \chi k \bigl(u(t_{m+\frac12})\vsigma(t_{m+\frac12}) - u^{m+\frac12}\vsigma^{m+\frac12},\nab e_u^{m+1}\bigr)\\\nonumber
			&:= III_1 + III_2.
		\end{align}
		
		Using Lemma \ref{Lemma_Holder} and  Lemma \ref{Stability_PDE} with the deterministic $(u_0,\sigma_0)$, we have
		\begin{align*}
			\mE[III_1] &=\chi \int_{t_m}^{t_{m+1}}\mE[\bigl(u(s)[\vsigma(s) - \vsigma(t_{m+\frac12})] + [u(s) - u(t_{m+\frac12})]\vsigma(t_{m+\frac12}), \nab e_{u}^{m+1}\bigr)]\, ds \\\nonumber
			&\leq \frac{C\chi^2}{\nu} \int_{t_m}^{t_{m+1}}\mE[\|u(s)[\vsigma(s) - \vsigma(t_{m+\frac12})]\|^2_{L^2} + \|[u(s) - u(t_{m+\frac12})]\vsigma(t_{m+\frac12})\|^2_{L^2}]\, ds \\\nonumber
			&\qquad\qquad\qquad+ \frac{\nu k}{64}\mE\bigl[\|\nab e_u^{m+1}\|^2_{L^2}\bigr]\\\nonumber
			&\leq CC_1k^2 + \frac{\nu k}{64}\mE\bigl[\|\nab e_u^{m+1}\|^2_{L^2}\bigr].
		\end{align*}
		
		To control $III_2$, we add and subtract $u^{m+\frac12}\vsigma(t_{m+\frac12})$ to get
		\begin{align*}
			III_2 &= \chi k \bigl(e_u^{m+\frac12}\vsigma(t_{m+\frac12}),\nab e_u^{m+1}\bigr) + \chi k \bigl(u^{m+\frac12}e_{\vsigma}^{m+\frac12},\nab e_u^{m+1}\bigr).
		\end{align*}	
		Using Ladyzhenskaya's inequality, we obtain
		\begin{align*}
			III_2 &\leq \frac{\nu k}{64} \|\nab e_u^{m+1}\|^2_{L^2} + \frac{C\chi^2k}{\nu}\bigl(\|e_u^{m+\frac12}\vsigma(t_{m+\frac12})\|^2_{L^2} + \|u^{m+\frac12}e_{\vsigma}^{m+\frac12}\|^2_{L^2}\bigr)\\\nonumber
			&\leq \frac{\nu k}{64} \|\nab e_u^{m+1}\|^2_{L^2} + \frac{C\chi^2k}{\nu}\bigl(\|e_u^{m+\frac12}\|^2_{L^4}\|\vsigma(t_{m+\frac12})\|^2_{L^4} + \|u^{m+\frac12}\|^2_{L^4}\|e_{\vsigma}^{m+\frac12}\|^2_{L^4}\bigr)\\\nonumber
			&\leq \frac{\nu k}{64} \|\nab e_u^{m+1}\|^2_{L^2} +\frac{C\chi^2 k}{\nu}C_L^2 \|e_u^{m+\frac12}\|_{L^2}\|\nab e_u^{m+\frac12}\|_{L^2}\|\vsigma(t_{m+\frac12})\|_{L^2}\|\nab \vsigma(t_{m+\frac12})\|_{L^2}\\\nonumber
			&\qquad+\frac{C\chi^2 k}{\nu} C_L^2\|u^{m+\frac12}\|_{L^2}\|\nab u^{m+\frac12}\|_{L^2}\|e_{\vsigma}^{m+\frac12}\|_{L^2}\|\nab e_{\vsigma}^{m+\frac12}\|_{L^2}\\\nonumber
			&\leq \frac{\nu k}{64} \|\nab e_u^{m+1}\|^2_{L^2} + \frac{\nu k}{64} \|\nab e_u^{m+\frac12}\|^2_{L^2} + \frac{C\chi^4C_L^4}{\nu^3}k\|e_u^{m+\frac12}\|^2_{L^2}\|\vsigma(t_{m+\frac12})\|^4_{H^1}\\\nonumber
			&\qquad+\frac{C\chi^2 k}{\nu} C_L^2\|u^{m+\frac12}\|^2_{H^1}\|e_{\vsigma}^{m+\frac12}\|^2_{H^1}.
		\end{align*}	
		Moreover, subtracting \eqref{eq_sigma} to \eqref{eq_semi_discrete_sigma}, we obtain
		\begin{align}\label{eq_3.31}
			&	\bigl(e_{\vsigma}^{m+1}, \vvarphi\bigr) + \bigl(\nab\cdot e_{\vsigma}^{m+1},\nab\cdot\vvarphi\bigr) + \bigl(rot\, e_{\vsigma}^{m+1},\, rot\, \vvarphi\bigr) \\\nonumber
			&= -\bigl(u(t_{m+1}) - u(t_m),\nab\cdot\vvarphi\bigr) - \bigl(e_{u}^m,\nab\cdot\vvarphi\bigr).
		\end{align}
		
		Therefore, for any $m\geq 1$, we have
		\begin{align}\label{eq_3.32}
			&	\bigl(e_{\vsigma}^{m+\frac12}, \vvarphi\bigr) + \bigl(\nab\cdot e_{\vsigma}^{m+\frac12},\nab\cdot\vvarphi\bigr) + \bigl(rot\, e_{\vsigma}^{m+\frac12},\, rot\, \vvarphi\bigr) \\\nonumber
			&= -\frac12\bigl(u(t_{m+1}) - u(t_m),\nab\cdot\vvarphi\bigr) \\\nonumber
			&\qquad - \frac12\bigl(u(t_{m}) - u(t_{m-1}),\nab\cdot\vvarphi\bigr) - \bigl(e_{u}^{m-\frac12},\nab\cdot\vvarphi\bigr).
		\end{align}
		
		Taking $\vvarphi = e_{\vsigma}^{m+\frac12}$ in \eqref{eq_3.32} we arrive at
		\begin{align*}
			\|e_{\vsigma}^{m+\frac12}\|^2_{H^1} \leq C\bigl[\|u(t_{m+1}) - u(t_m)\|^2_{L^2} +\|u(t_m) - u(t_{m-1})\|^2_{L^2}\bigr] + C\|e_u^{m-\frac12}\|^2_{L^2}.
		\end{align*}
		With this, we update the estimates for $III_2$ as follows:
		\begin{align*}
			III_2 	&\leq \frac{\nu k}{64} \|\nab e_u^{m+1}\|^2_{L^2} + \frac{\nu k}{64} \|\nab e_u^{m+\frac12}\|^2_{L^2} + \frac{C\chi^4C_L^4}{\nu^3}k\|e_u^{m+\frac12}\|^2_{L^2}\|\vsigma(t_{m+\frac12})\|^4_{H^1}\\\nonumber
			&\qquad+\frac{C\chi^2 k}{\nu} C_L^2\|u^{m+\frac12}\|^2_{H^1}\|e_u^{m-\frac12}\|^2_{L^2} \\\nonumber
			&\qquad\qquad+ \frac{C\chi^2 k}{\nu} C_L^2\|u^{m+\frac12}\|^2_{H^1}\|u(t_{m+1}) - u(t_m)\|^2_{L^2} \\\nonumber
			&\qquad\qquad+ \frac{C\chi^2 k}{\nu} C_L^2\|u^{m+\frac12}\|^2_{H^1}\|u(t_{m}) - u(t_{m-1})\|^2_{L^2}.
		\end{align*}
		
		With this and using Lemma \ref{Lemma_Holder}, Lemma \ref{Lemma_stability}, and Lemma \ref{Stability_PDE}, we obtain
		\begin{align}
			\mE[III_2] &\leq  \frac{\nu k}{32} \mE[\|\nab e_u^{m+1}\|^2_{L^2} + \frac{CC_1\chi^4C_L^4}{\nu^3}k\mE[\|e_u^{m+\frac12}\|^2_{L^2}] \\\nonumber
			&\qquad+\frac{CC_1\chi^2 k}{\nu} C_L^2\mE[\|e_u^{m-\frac12}\|^2_{L^2}] + \frac{CC_1\chi^2C_L^2K_2}{\nu}k^2.
		\end{align}
		
		Next, we analyze the noise term as follows. 
		\begin{align}
			IV&= \delta\Bigl(\int_{t_m}^{t_{m+1}} (\vb\cdot\nab u(s) - \vb\cdot\nab u^m)\, dW(s), e^{m+1}_u\Bigr) \\\nonumber
			&\qquad-\frac{\delta}{2}\bigl((\vb\cdot\nab u^{m+1} - \vb\cdot\nab u^m)\Delta W_m, e_u^{m+1}\bigr)\\\nonumber
			&\qquad+ \frac{\delta^2}{2}\int_{t_m}^{t_{m+1}}\bigl(\nab\cdot(B\nab u(s)),e_u^{m+1}\bigr)\, ds\\\nonumber
			&:= IV_1 + IV_2 + IV_3.
		\end{align}
		
		Using the martingale property of It\^o integral, the It\^o isometry and Lemma \ref{Lemma_Holder}, we have
		\begin{align}
			\mE[IV_1] &= \delta\mE\Bigl[\Bigl(\int_{t_m}^{t_{m+1}} (\vb\cdot\nab u(s) - \vb\cdot\nab u^m)\, dW(s), e^{m+1}_u - e_u^m\Bigr) \Bigr]\\\nonumber
			&= \delta\mE\Bigl[\Bigl(\int_{t_m}^{t_{m+1}} (\vb\cdot\nab u(s) - \vb\cdot\nab u(t_m))\, dW(s), e^{m+1}_u - e_u^m\Bigr) \Bigr]\\\nonumber
			&\qquad+ \delta\mE\Bigl[\Bigl( \vb\cdot\nab e_u^{m}, e^{m+1}_u - e_u^m\Bigr)\Delta W_m \Bigr]\\\nonumber
			&\leq Ck^2 + \frac{1}{64}\mE\bigl[\|e_u^{m+1} - e_u^m\|^2_{L^2}\bigr] \\\nonumber
			&\qquad+ \delta^2k\mE\bigl[\|\vb\cdot\nab e^m_u\|^2_{L^2}\bigr] + \frac14\mE\bigl[\|e_u^{m+1} - e_u^m\|^2_{L^2}\bigr].
		\end{align}
		
		For $IV_2$ and $IV_3$, we use the integration by parts and then the scheme \eqref{eq_semi_discrete_u} to control them as follow: 
		\begin{align*}
			IV_2 + IV_3 &= \frac{\delta}{2}\bigl(u^{m+1} - u^m,\vb\cdot\nab e_u^{m+1}\bigr)\Delta W_m+ \frac{\delta^2}{2}\int_{t_m}^{t_{m+1}}\bigl(\nab\cdot(B\nab u(s)),e_u^{m+1}\bigr)\, ds\\\nonumber
			&=\frac{\delta}{2}\Bigl[-\nu k\bigl(\nab u^{m+\frac12},\nab(\vb\cdot\nab e_u^{m+1})\bigr)  +\chi k \bigl(u^{m+\frac12}\vsigma^{m+\frac12}, \nab(\vb\cdot\nab e_u^{m+1})\bigr)  \\\nonumber
			&\qquad+ \delta \bigl(\vb\cdot\nab u^{m+\frac12}, \vb\cdot\nab e_u^{m+1}\bigr)\Delta W_m\Bigr]\Delta W_m+ \frac{\delta^2}{2}\int_{t_m}^{t_{m+1}}\bigl(\nab\cdot(B\nab u(s)),e_u^{m+1}\bigr)\, ds\\\nonumber
			&=-\frac{\delta}{2}\nu k\bigl(\nab u^{m+\frac12},\nab(\vb\cdot\nab e_u^{m+1})\bigr) \Delta W_m +\frac{\delta}{2}\chi k \bigl(u^{m+\frac12}\vsigma^{m+\frac12}, \nab(\vb\cdot\nab e_u^{m+1})\bigr) \Delta W_m \\\nonumber
			&\qquad+ \frac{\delta^2}{2} \bigl(\vb\cdot\nab u^{m+\frac12}, \vb\cdot\nab e_u^{m+1}\bigr) (\Delta W_m)^2+ \frac{\delta^2}{2}\int_{t_m}^{t_{m+1}}\bigl(\nab\cdot(B\nab u(s)),e_u^{m+1}\bigr)\, ds\\\nonumber
			&=-\frac{\delta}{2}\nu k\bigl(\nab u^{m+\frac12},\nab(\vb\cdot\nab e_u^{m+1})\bigr) \Delta W_m +\frac{\delta}{2}\chi k \bigl(u^{m+\frac12}\vsigma^{m+\frac12}, \nab(\vb\cdot\nab e_u^{m+1})\bigr) \Delta W_m \\\nonumber
			&\qquad-\frac{\delta^2}{2} \bigl(\nab\cdot(B\nab u^{m+\frac12}),  e_u^{m+1}\bigr) (\Delta W_m)^2 + \frac{\delta^2}{2}\int_{t_m}^{t_{m+1}}\bigl(\nab\cdot(B\nab u(s)),e_u^{m+1}\bigr)\, ds\\\nonumber
			&=-\frac{\delta}{2}\nu k\bigl(\nab u^{m+\frac12},\nab(\vb\cdot\nab e_u^{m+1})\bigr) \Delta W_m +\frac{\delta}{2}\chi k \bigl(u^{m+\frac12}\vsigma^{m+\frac12}, \nab(\vb\cdot\nab e_u^{m+1})\bigr) \Delta W_m \\\nonumber
			&\qquad +\frac{\delta^2}{2} \bigl(\nab\cdot(B\nab u^{m+\frac12}),  e_u^{m+1}\bigr) \Bigl[k-(\Delta W_m)^2\Bigr] \\\nonumber
			&\qquad+ \frac{\delta^2}{2}\int_{t_m}^{t_{m+1}}\bigl(\nab\cdot(B\nab u(s) - B\nab u(t_{m+\frac12})),e_u^{m+1}\bigr)\, ds \\\nonumber
			&\qquad+ \frac{\delta^2}{2}k\bigl(\nab\cdot(B\nab e_u^{m+\frac12}),e_u^{m+1}\bigr)\\\nonumber
			&:= IV_{2,3}^1 +  IV_{2,3}^2+ IV_{2,3}^3+ IV_{2,3}^4 +IV_{2,3}^5.
		\end{align*}
		
		Using integration by parts, we obtain
		\begin{align*}
			IV_{2,3}^1 &= \frac{\delta}{2} \nu k\bigl(\Delta u^{m+\frac12}, \vb\cdot\nab e_u^{m+1}\bigr)\Delta W_m\\\nonumber
			&\leq \frac{\nu k}{64} \|\nab e_u^{m+1}\|^2_{L^2} + C|\vb|^2\delta^2 k\|\Delta u^{m+\frac12}\|^2_{L^2}|\Delta W_m|^2\\\nonumber
			&\leq \frac{\nu k}{64} \|\nab e_u^{m+1}\|^2_{L^2} + C|\vb|^2\delta^2 k^2\|\Delta u^{m+\frac12}\|^4_{L^2} + C|\vb|^2\delta^2|\Delta W_m|^4,
		\end{align*}
		which, together with Lemma \ref{Lemma_stability}, implies that
		\begin{align*}
			\mE[IV_{2,3}^1] &\leq \frac{\nu k}{64} \mE\bigl[\|\nab e_u^{m+1}\|^2_{L^2}\bigr] + C|\vb|^2\delta^2 k^2\mE\bigl[\|\Delta u^{m+\frac12}\|^4_{L^2}\bigr] + C|\vb|^2\delta^2k^2\\\nonumber
			&\leq \frac{\nu k}{64} \mE\bigl[\|\nab e_u^{m+1}\|^2_{L^2}\bigr] + C|\vb|^2\delta^2 (1+K_3^2)k^2.
		\end{align*}
		
		Similarly, using integration by parts and the Ladyzhenskaya inequality, we also obtain
		\begin{align*}
			IV_{2,3}^2 &= -\frac{\delta}{2}\chi k \bigl(\nab\cdot(u^{m+\frac12}\vsigma^{m+\frac12}), \vb\cdot \nab e_u^{m+1}\bigr)\Delta W_{m}\\\nonumber
			&\leq \frac{\nu k}{64} \|\nab e_u^{m+1}\|^2_{L^2} + \frac{C\chi^2\delta^2}{\nu}k\|\nab\cdot(u^{m+\frac12}\vsigma^{m+\frac12})\|^2_{L^2}|\Delta W_m|^2\\\nonumber
			&= \frac{\nu k}{64} \|\nab e_u^{m+1}\|^2_{L^2} + \frac{C\chi^2\delta^2}{\nu}k\|\nab u^{m+\frac12}\vsigma^{m+\frac12} + u^{m+\frac12}\nab\cdot \vsigma^{m+\frac12}\|^2_{L^2}|\Delta W_m|^2\\\nonumber
			&\leq	\frac{\nu k}{64} \|\nab e_u^{m+1}\|^2_{L^2} + \frac{C\chi^2\delta^2}{\nu}k\Bigl[\|\nab u^{m+\frac12}\vsigma^{m+\frac12}\|^2_{L^2} \\\nonumber
			&\qquad\qquad\qquad\qquad+ \|u^{m+\frac12}\nab\cdot \vsigma^{m+\frac12}\|^2_{L^2}\Bigr]|\Delta W_m|^2\\\nonumber
			&\leq	\frac{\nu k}{64} \|\nab e_u^{m+1}\|^2_{L^2} + \frac{C\chi^2\delta^2}{\nu}k\Bigl[\|\nab u^{m+\frac12}\|^2_{L^4}\|\vsigma^{m+\frac12}\|^2_{L^4} \\\nonumber
			&\qquad\qquad\qquad\qquad+ \|u^{m+\frac12}\|^2_{L^4}\|\nab\cdot \vsigma^{m+\frac12}\|^2_{L^4}\Bigr]|\Delta W_m|^2\\\nonumber
			&\leq	\frac{\nu k}{64} \|\nab e_u^{m+1}\|^2_{L^2} + \frac{C\chi^2\delta^2}{\nu}C^2_Lk\Bigl[\| u^{m+\frac12}\|^2_{H^2}\|\vsigma^{m+\frac12}\|^2_{H^1} \\\nonumber
			&\qquad\qquad\qquad\qquad+ \|u^{m+\frac12}\|^2_{H^1}\|\nab\cdot \vsigma^{m+\frac12}\|_{L^2}\|\nab(\nab\cdot\vsigma^{m+\frac12})\|_{L^2}\Bigr]|\Delta W_m|^2\\\nonumber
			&\leq	\frac{\nu k}{64} \|\nab e_u^{m+1}\|^2_{L^2} + \frac{C\chi^2\delta^2}{\nu}C^2_Lk^2\Bigl[K_3^2K_1^2 + K_2^2K_1K_2\Bigr] + \frac{C\chi^2\delta^2}{\nu}C^2_L|\Delta W_m|^4,
		\end{align*}
		which implies
		\begin{align*}
			\mE[IV^2_{2,3}] &\leq 	\frac{\nu k}{64} \mE\bigl[\|\nab e_u^{m+1}\|^2_{L^2}\bigr] + \frac{C\chi^2\delta^2}{\nu}C^2_Lk^2\bigl[K_1^2K_3^2 + K_1K_2^3\bigr] +  \frac{C\chi^2\delta^2}{\nu}C^2_L k^2.
		\end{align*}
		
		To estimate $IV_{2,3}^3$, we use the independence of $\Delta W_m$ and $e_u^m, u^m$ to obtain
		\begin{align*}
			\mE[IV_{2,3}^3] &= \frac{\delta^2}{2}\mE\Bigl[ \bigl(\nab\cdot(B\nab u^{m+\frac12}),  e_u^{m+1}  - e_u^m\bigr) \bigl[k-(\Delta W_m)^2\bigr] \Bigr] \\\nonumber
			&+ \frac{\delta^2}{4}\mE\Bigl[ \bigl(\nab\cdot(B\nab u^{m+1}- B\nab u^{m}), e_u^m\bigr) \bigl[k-(\Delta W_m)^2\bigr] \Bigr]\\\nonumber
			&:= IV_{2,3}^{3a} + IV_{2,3}^{3b}.
		\end{align*}
		Using Lemma \ref{Lemma_stability}, the term $IV_{2,3}^{3a}$ can be estimated as
		\begin{align*}
			IV_{2,3}^{3a} &\leq \frac{1}{64}\mE\bigl[\|e_u^{m+1} - e^m_u\|^2_{L^2}\bigr] + Ck^2 \mE\bigl[\|\nab\cdot(B\nab u^{m+\frac12})\|^4_{L^2}\bigr] + C\frac{1}{k^2}\mE\bigl[|(\Delta W_m)^2 - k|^4\bigr]\\\nonumber
			&\leq  \frac{1}{64}\mE\bigl[\|e_u^{m+1} - e^m_u\|^2_{L^2}\bigr] + Ck^2 \mE\bigl[\|\nab\cdot(B\nab u^{m+\frac12})\|^4_{L^2}\bigr] + Ck^2\\\nonumber
			&\leq  \frac{1}{64}\mE\bigl[\|e_u^{m+1} - e^m_u\|^2_{L^2}\bigr] + C(1+K_3^2)k^2 
		\end{align*}
		
		Using integration by parts and Lemma \ref{Lemma_Holder}(c), we obtain
		\begin{align*}
			IV_{2,3}^{3b} &= \frac{\delta^2}{4}\mE\bigl[\bigl(u^{m+1}  - u^m, \nab\cdot (B\nab e_u^m)\bigr)[k-(\Delta W_m)^2]\bigr]\\\nonumber
			&= -\frac{\delta^2}{4}\mE\bigl[\bigl(e_{u}^{m+1}  - e_u^m,\nab\cdot(B\nab e_u^m)\bigr)[k-(\Delta W_m)^2]\bigr] \\\nonumber&\qquad+ \frac{\delta^2}{4}\mE\bigl[\bigl({u(t_{m+1})}  - u(t_{m}), \nab\cdot(B\nab e_u^m)\bigr)[k-(\Delta W_m)^2]\bigr]\\\nonumber
			&= -\frac{\delta^2}{4}\mE\bigl[\bigl(e_{u}^{m+1}  - e_u^m,\nab\cdot(B\nab e_u^m)\bigr)[k-(\Delta W_m)^2]\bigr] \\\nonumber&\qquad- \frac{\delta^2}{4}\mE\bigl[\bigl({B\nab u(t_{m+1})}  - B\nab u(t_{m}), \nab e_u^m\bigr)[k-(\Delta W_m)^2]\bigr]\\\nonumber
			&\leq \frac{1}{64}\mE\bigl[\|e_u^{m+1} - e_u^{m}\|^2_{L^2}\bigr] + \frac{\nu k}{64}\mE\bigl[\|\nab e_u^m\|^2_{L^2}\bigr] \\\nonumber
			&\qquad+ C\mE\bigl[\|\nab\cdot(B\nab e_u^{m})\|^2_{L^2}[k - (\Delta W_m)^2]^2\bigr] \\\nonumber
			&\qquad+ \frac{C}{k}\mE\bigl[\|\nab(u(t_{m+1}) - u(t_m))\|^2_{L^2}[k-(\Delta W_m)^2]^2\bigr]\\\nonumber
			&\leq \frac{1}{64}\mE\bigl[\|e_u^{m+1} - e_u^{m}\|^2_{L^2}\bigr] + \frac{\nu k}{64}\mE\bigl[\|\nab e_u^m\|^2_{L^2}\bigr] \\\nonumber
			&\qquad+ C(C_2+K_3)k^2 + C\mE\bigl[\|\nab(u(t_{m+1}) - u(t_m))\|^4_{L^2}\bigr] + \frac{C}{k^2}\mE\bigl[k-(\Delta W_m)^2]^4\bigr]\\\nonumber
			&\leq \frac{1}{64}\mE\bigl[\|e_u^{m+1} - e_u^{m}\|^2_{L^2}\bigr] + \frac{\nu k}{64}\mE\bigl[\|\nab e_u^m\|^2_{L^2}\bigr] + C(1+C_2+K_3)k^2.
		\end{align*}
		
		Next, to estimate $IV_{2,3}^4$, we use integration by parts and Lemma \ref{Lemma_Holder} to get
		\begin{align*}
			\mE[IV_{2,3}^4] &= -\frac{\delta^2}{2}\int_{t_m}^{t_{m+1}}\mE\bigl[\bigl(B\nab(u(s) - u(t_{m+\frac12})), \nab e_u^{m+1}\bigr)\bigr]\, ds\\\nonumber
			&\leq Ck^2 + \frac{\nu k}{64}\mE\bigl[\|\nab e_u^{m+1}\|^2_{L^2}\bigr].
		\end{align*}
		
		Using integration by parts, Lemma \ref{Stability_PDE}, and Lemma \ref{Lemma_stability}, we get
		\begin{align*}
			\mE[IV_{2,3}^5] &= -\frac{\delta^2}{2}k\mE[\bigl(\vb\cdot\nab e_u^{m+\frac12},\vb\cdot\nab e_u^{m+1}\bigr)]\\\nonumber
			&=-\frac{\delta^2}{4}k\mE[\|\vb\cdot\nab e_u^{m+1}\|^2_{L^2}] - \frac{\delta^2}{4}k\mE[\|\vb\cdot\nab e_u^m\|^2_{L^2}] \\\nonumber&\qquad+\frac{\delta^2}{4}k\mE[\bigl(\nab\cdot(B\nab e_u^m), e_u^{m+1} - e_u^m\bigr)]\\\nonumber
			&\leq -\frac{\delta^2}{4}k\mE[\|\vb\cdot\nab e_u^{m+1}\|^2_{L^2}] - \frac{\delta^2}{4}k\mE[\|\vb\cdot\nab e_u^m\|^2_{L^2}] \\\nonumber&\qquad+\frac{1}{64}\mE[\|e_u^{m+1} - e_u^m\|^2_{L^2}] + Ck^2\mE[\|\nab\cdot(B\nab e_u^m)\|^2_{L^2}]\\\nonumber
			&\leq -\frac{\delta^2}{4}k\mE[\|\vb\cdot\nab e_u^{m+1}\|^2_{L^2}] - \frac{\delta^2}{4}k\mE[\|\vb\cdot\nab e_u^m\|^2_{L^2}] \\\nonumber&\qquad+\frac{1}{64}\mE[\|e_u^{m+1} - e_u^m\|^2_{L^2}] + C(C_2+K_3)k^2.
		\end{align*}
		
		Now, substituting all the estimates into \eqref{eq_3.29} with the expectation, then rearrange terms, we obtain
		\begin{align}\label{eq_3.35}
			&\frac12\mE\bigl[\|e^{m+1}_u\|^2_{L^2} - \|e_u^m\|^2_{L^2}\bigr] +\frac{5}{32}\mE[\|e_u^{m+1} - e_u^m\|^2_{L^2}]  \\\nonumber
			&\qquad\qquad+\frac{11\nu k}{64}\mE[\|\nab e_u^{m+1}\|^2_{L^2}] +\left( \frac{15\nu }{32} - \frac{\delta^2 |\vb|^2}{2}\right)k\mE[\|\nab e_u^m\|^2_{L^2}]\\\nonumber
			&\qquad\qquad+ \frac{\delta^2 k}{4}\mE[\|\vb\cdot\nab e_u^{m+1}\|^2_{L^2}] -\frac{\delta^2 k}{4}\mE[\|\vb\cdot\nab e_u^{m}\|^2_{L^2}]
			\\\nonumber
			&\leq   \frac{CC_1\chi^4C_L^4}{\nu^3}k\mE[\|e_u^{m+\frac12}\|^2_{L^2}] +\frac{CC_1\chi^2 C_L^2}{\nu}k \mE[\|e_u^{m-\frac12}\|^2_{L^2}]+ CKk^2,
		\end{align}
		where $\frac{15\nu}{32} - \frac{\delta^2|\vb|}{2}>0$ by using the assumption \eqref{assume_delta} of $\delta$, and 
		\begin{align*}
			K: = &\left\{1+C_1+C_2+K_3 + \frac{C_1\chi^2C_L^2K_2}{\nu} + |\vb|^2\delta^2(1+K_3^2)\right.\\\nonumber
			&\left.+\frac{\chi^2C_L^2\delta^2}{\nu}[K_3^2K_1^2 + K_1K_2^3] + \frac{\chi^2\delta^2}{\nu}+ K_3^2\right\}.
		\end{align*}
		
		Next, applying the summation $\sum_{m=1}^{\ell}$ for any $1\leq \ell \leq M-1$ to \eqref{eq_3.35} we obtain
		\begin{align}\label{eq_3.36}
			&\frac12\mE\bigl[\|e^{\ell+1}_u\|^2_{L^2}\bigr] +\frac{5}{32}\sum_{m=1}^{\ell}\mE[\|e_u^{m+1} - e_u^m\|^2_{L^2}]  + \frac{\delta^2}{4}k\mE[\|\vb\cdot\nab e_u^{\ell+1}\|^2_{L^2}]\\\nonumber
			&\qquad+\frac{11\nu }{64}k\sum_{m=1}^{\ell}\mE[\|\nab e_u^{m+1}\|^2_{L^2}] + \Bigl(\frac{15\nu}{64} - \frac{\delta^2|\vb|}{2}\Bigr)k\sum_{m=1}^{\ell}\mE[\|\nab e_u^m\|^2_{L^2}]\\\nonumber
			&\leq \frac12\mE[\|e_u^1\|^2_{L^2}]  +  \frac{\delta^2k}{4}\mE[\|\vb\cdot\nab e_u^1\|^2_{L^2}] \\\nonumber
			&\qquad+CKTk + \frac{CC_1\chi^2C_L^2}{\nu}\left(1+\frac{\chi^2C_L^2}{\nu^2}\right)k\sum_{m=0}^{\ell}\mE[\|e_u^{m+1}\|^2_{L^2} + \|e_u^{m}\|^2_{L^2}]\\\nonumber
			&\leq \frac12\mE[\|e_u^1\|^2_{L^2}]  +  C(KT + \delta^2|b|^2(K_1+C_1))k \\\nonumber
			&\qquad + \frac{CC_1\chi^2C_L^2}{\nu}\left(1+\frac{\chi^2C_L^2}{\nu^2}\right)k\sum_{m=0}^{\ell}\mE[\|e_u^{m+1}\|^2_{L^2} + \|e_u^{m}\|^2_{L^2}],
		\end{align}
		where the last inequality of \eqref{eq_3.36} is obtained by using $\mE[\|\vb\cdot\nab e^1_u\|^2_{L^2}] \leq |b|^2(K_1 +C_1)$.
		
		Applying the standard discrete Gronwall inequality to \eqref{eq_3.36}, we arrive at
		\begin{align}\label{3.37}
			&\frac12\mE\bigl[\|e^{\ell+1}_u\|^2_{L^2}\bigr] +\frac{5}{32}\sum_{m=1}^{\ell}\mE[\|e_u^{m+1} - e_u^m\|^2_{L^2}]  + \frac{\delta^2}{4}k\mE[\|\vb\cdot\nab e_u^{\ell+1}\|^2_{L^2}]\\\nonumber
			&\qquad+\frac{11\nu }{64}k\sum_{m=1}^{\ell}\mE[\|\nab e_u^{m+1}\|^2_{L^2}] + \Bigl(\frac{15\nu}{64} - \frac{\delta^2|\vb|}{2}\Bigr)k\sum_{m=1}^{\ell}\mE[\|\nab e_u^m\|^2_{L^2}]\\\nonumber
			&\leq \left(\frac12\mE[\|e_u^1\|^2_{L^2}] + C(KT + \delta^2|b|^2(K_1+C_1))k\right)\exp\left[\frac{CTC_1\chi^2C_L^2}{\nu}\left(1+\frac{\chi^2C_L^2}{\nu^2}\right)\right].
		\end{align}
		
		It is left to estimate $\mE[\|e_u^1\|^2_{L^2}]$. To do that, 
		Taking $m = 0$ in \eqref{eq_3.28} and use the fact that $e_u^0 =0$, and $e_{\vsigma}^0=0$, we get
		\begin{align}\label{eq3.39}
			\|e_u^1\|^2_{L^2} +\frac{\nu k}{2}\|\nab e_u^1\|^2_{L^2} &= \int_{t_0}^{t_1} \bigl(\nab(u(t_{1/2}) - u(s)), \nab e_u^1\bigr)\, ds \\\nonumber
			&+\chi  \int_{t_0}^{t_1} \bigl(u(s)\vsigma(s) - u^{\frac12}\vsigma^{\frac12}, \nab e^1_u\bigr)\, ds\\\nonumber
			&+ \delta \Bigl(\int_{t_0}^{t_{1}} (\vb\cdot\nab u(s) - \vb\cdot\nab u^{\frac12})\circ dW(s), e_u^1\Bigr)\\\nonumber
			&:= L_1 + L_2 + L_3.
		\end{align}
		
		Using Lemma \ref{Lemma_Holder}, we have
		\begin{align*}
			\mE[L_1] &\leq Ck^2 + \frac{\nu k}{8}\mE[\|\nab e_u^1\|^2_{L^2}].
		\end{align*}
		
		Using Lemma \ref{Stability_PDE}, Lemma \ref{Lemma_stability} and Lemma \ref{Lemma_Holder},
		\begin{align*}
			\mE[L_2]&= \chi  \int_{t_0}^{t_1}\mE[\bigl(u(s)\vsigma(s) - u(t_{\frac12})\vsigma(t_{\frac12}), \nab e^1_u\bigr)]\, ds +  \frac{\chi k}{2} \mE[\bigl(u(t_{\frac12})\vsigma(t_{\frac12}) - u^{\frac12}\vsigma^{\frac12}, \nab e^1_u\bigr)] \\\nonumber
			&\leq \frac{\nu k}{8} \mE[\|\nab e_u^1\|^2_{L^2}] + Ck^2 + C(C_1+ K_2)k.
		\end{align*}
		
		To estimate $L_3$, we proceed as below. Using Lemma \ref{Stability_PDE} and Lemma \ref{Lemma_Holder}, we obtain
		\begin{align*}
			\mE[L_3]& = \delta \mE\Bigl[\bigl(\vb\cdot\nab e_u^{\frac12}, e_u^1\bigr)\Delta W_0+ \delta \Bigl(\int_{t_0}^{t_{1}} (\vb\cdot\nab u(s) - \vb\cdot\nab u(t_{\frac12}))\circ dW(s), e_u^1\Bigr)\Bigr]\\\nonumber
			& = 0+\delta \mE\Bigl[\Bigl(\int_{t_0}^{t_{1}} (\vb\cdot\nab u(s) - \vb\cdot\nab u(t_{\frac12}))\circ dW(s), e_u^1\Bigr)\Bigr]\\\nonumber
			& =\delta \mE\Bigl[\Bigl(\int_{t_0}^{t_{1}} (\vb\cdot\nab u(s) - \vb\cdot\nab u(t_{\frac12})) dW(s), e_u^1\Bigr)\Bigr] \\\nonumber
			&\qquad+ \frac{\delta^2}{2}\mE\Bigl[ \int_{t_0}^{t_{1}} \bigl(\nab\cdot(B\nab u(s)), e_u^1\bigr)\, ds\Bigr]\\\nonumber
			& =\delta \mE\Bigl[\Bigl(\int_{t_0}^{t_{1}} (\vb\cdot\nab u(s) - \vb\cdot\nab u(t_0)) dW(s), e_u^1\Bigr)\Bigr] \\\nonumber
			&\qquad-\frac{\delta}{2} \mE\bigl[\bigl( (\vb\cdot\nab u(t_1) - \vb\cdot\nab u(t_0)) \Delta W_0, e_u^1\bigr)\bigr] \\\nonumber
			&\qquad+ \frac{\delta^2}{2}\mE\Bigl[ \int_{t_0}^{t_{1}} \bigl(\nab\cdot(B\nab u(s)), e_u^1\bigr)\, ds\Bigr]\\\nonumber
			&\leq Ck^2 + \frac{1}{4}\mE[\|e_u^1\|^2_{L^2}]  + C\mE[\|\vb\cdot\nab(u(t_1) - u(t_0))\|^2_{L^2}|\Delta W_0|^2]\\\nonumber
			&\qquad+ Ck \mE\Bigl[\int_{t_0}^{t_1} \|\nab\cdot(B\nab u(s))\|^2_{L^2}\, ds\Bigr] + \frac{\nu k}{8}\mE[\|\nab e_u^1\|^2_{L^2}]\\\nonumber
			&\leq C(1+C_2)k+ \frac{1}{4}\mE[\|e_u^1\|^2_{L^2}]  + \frac{\nu k}{8}\mE[\|\nab e_u^1\|^2_{L^2}].
		\end{align*}
		
		Substituting these estimates back into \eqref{eq3.39}, we obtain
		\begin{align}\label{eq_3.40}
			\mE[\|e_u^1\|^2_{L^2}] +\frac{\nu k}{2}\mE[\|\nab e_u^1\|^2_{L^2}] &\leq C(1+C_1+C_2+K_2)k.
		\end{align}
		Finally, we obtain \eqref{estimate_u} by combining \eqref{eq_3.36} and \eqref{eq_3.40}. The proof is complete.
		
	\end{proof}

	\section{Fully discrete splitting mixed finite element method}\label{sec4}
	
	\subsection{Formulation of the finite element method} In this part, we consider a semi-implicit fully discrete finite element method for Algorithm 1.
	
	Let $\{\mathcal{T}_h\}_{h>0}$ be a quasi-uniform triangular mesh of $D\subset \mathbb{R}^2$, with the meshsize $0<h\ll 1$. We consider the finite element spaces $X^h_u\times X^h_{\vsigma}\times X^h_{c}$, which are generated by polynomials of degree at most $r_1, r_2,r_3\geq 1$. In addition, the splitting mixed element method doesn't require the LBB Inf-Sup condition for its finite element spaces. So, we freely choose the linear polynomials for our finite element spaces. Although we consider only linear finite element spaces in this paper, all the error estimates from the finite element approximation can be easily extended to higher polynomial finite element spaces.
	
	Next, we introduce the following $L^2$ projection: 
	\begin{align}
		&\cP_u: L_{per}^2(D) \rightarrow X_u^h, \qquad \cP_{\vsigma}: \vH_{per}^1(D) \rightarrow X_{\vsigma}^h, \\\nonumber
		&\cP_c: L_{per}^2(D) \rightarrow X_c^h,
	\end{align}
	where the operators satisfy
	\begin{align*}
		\bigl(v - \cP_{u} v, \varphi_h\bigr) &= 0\qquad\forall \varphi_h \in X_u^h,\\\nonumber
		\bigl(\vv - \cP_{\vsigma} \vv, \vvarphi_h\bigr) + \bigl(\nab\cdot(\vv - \cP_{\vsigma} \vv), \nab\cdot\vvarphi_h\bigr) +\bigl(rot(\vv - \cP_{\vsigma} \vv), rot\vvarphi_h\bigr) &= 0\qquad\forall \vvarphi_h \in X_{\vsigma}^h,\\\nonumber
		\bigl(v - \cP_{c} v, \varphi_h\bigr) &= 0\qquad\forall \varphi_h \in X_c^h.
	\end{align*}

	Then, we have the following interpolation inequalities:
	\begin{align}\label{ineq_interpolation}
		\|v- \cP_{u} v\|_{L^2} + h\|v - \cP_{u}v\|_{H^1} &\leq Ch^{2}\|v\|_{H^{2}}\qquad \forall v \in H_{per}^2(D),\\\nonumber
		\|\pphi- \cP_{\vsigma} \pphi\|_{L^2} + h\|\pphi - \cP_{\vsigma}\pphi\|_{\vH^1} &\leq Ch^{2}\|\pphi\|_{\vH^{2}}\qquad\forall \pphi\in \vH_{per}^2(D),\\\nonumber
		\|v- \cP_{c} v\|_{L^2} + h\|v - \cP_{c}v\|_{H^1} &\leq Ch^{2}\|v\|_{H^{2}}\qquad\forall v \in H^2_{per}(D).
	\end{align}

	Then we can state the fully discrete Crank-Nicolson finite element method for \eqref{Weak_formulation} in the following algorithm.
	
	\bigskip
	{\bf Algorithm 2:} Let $u^0_h = \mathcal{P}_u u_0$ and $\vsigma_h^0 = \cP_{\vsigma} \vsigma_0$. For any $0\leq m \leq M-1$, we do following iterations:
	Find $\bigl(u^{m+1}_h, \vsigma_h^{m+1}, c_h^{m+1}\bigr) \in X^h_u\times X^h_{\vsigma}\times X^h_{c}$ such that $\mP$- a.s.
	\begin{align}
		\label{eq_discrete_sigma}	&\bigl(\vsigma_h^{m+1},\vvarphi_h\bigr) + \bigl(\nab\cdot\vsigma_h^{m+1},\nab\cdot\vvarphi_h\bigr) + \bigl(rot\, \vsigma_h^{m+1},rot\, \vvarphi_h\bigr) = -\bigl(u_h^{m}, \nab\cdot\vvarphi_h\bigr),\\
		\medskip
		\label{eq_discrete_u}	&\bigl(u_h^{m+1}- u_h^{m},v_h\bigr) + k \nu\bigl(\nab u_h^{m+\frac12},\nab v_h\bigr) = k \chi\bigl(u_h^{m+\frac12}\vsigma_h^{m+\frac12},\nab v_h\bigr) \\\nonumber
		&\qquad\qquad\qquad\qquad\qquad\qquad\qquad\qquad\qquad+ \delta\bigl(\vb\cdot\nab u_h^{m+\frac12}\Delta W_m, v_h\bigr),\\\medskip
		\label{eq_discrete_c}	&\bigl(c_h^{m+1}, \psi_h\bigr) = \bigl(\nab\cdot\vsigma_h^{m+1},\psi_h\bigr) + \bigl(u^{m+1}_h,\psi_h\bigr),
	\end{align}
	for all $\bigl(v_h, \vvarphi_h,\psi_h\bigr) \in X^h_u\times X^h_{\vsigma}\times X^h_{c}$.
	
	\medskip
	
	Since we choose $u^0_h = \mathcal{P}_u u_0$ and $\vsigma_h^0 = \cP_{\vsigma} \vsigma_0$ and \eqref{ineq_interpolation}, without loss of generality, we assume that $u^0_h = u_0$ and $\vsigma_h^0  = \vsigma_0$ in all the proofs of the stability and error estimates of Algorithm 2. 
	
	First, we can easily establish the mass-conservation for $\{u_h^m\}$. 
	\begin{lemma} $\{u^m_h\}$ satisfies the mass-conservation:
		\begin{align*}
			\int_D u^{m+1}_h = \int_D u_h^m= \cdots =\int_D u_h^0.
		\end{align*}
	\end{lemma}
	
	Next, we present the stability estimates of Algorithm 2. Its proof is similar to the proof of Lemma \ref{Lemma_stability} (a). So, we skip it to save space.
	
	\begin{lemma}\label{lemma_fem_stability} Let $(u^m_h,\vsigma_h^m,c^m_h)$ be the solution of Algorithm 2. Assume that $u_0 \in L^2_{per}$ satisifes the condition \eqref{stability_criteria}. There exist a constant $\tilde{K} = \tilde{K}(u_0)>0$ such that $\mP$-a.s.
		\begin{enumerate}
			\item[(a)] $\displaystyle \max_{1\leq m \leq M}\|u_h^m\|^2_{L^2} + \nu k\sum_{m=1}^{M-1}\|\nab u_h^{m+\frac12}\|^2_{L^2} \leq \tilde{K}$,
			\item[(b)] $\displaystyle \max_{1\leq m \leq M}\|\vsigma_h^m\|^2_{H^1} \leq \tilde{K}$,
			\item[(c)] $\displaystyle \max_{1\leq m \leq M} \|c_h^m\|^2_{L^2}  \leq \tilde{K}$.
		\end{enumerate}
	\end{lemma}
	
	Next, we state and prove the error estimates of Algorithm 2. 
	
	\subsection{Error estimates} In this part, we analyze and derive the error estimates for Algorithm 2.
	\begin{theorem}\label{Thm_error} Let $(u^m, \vsigma^m, c^m)$ and $(u_h^m, \vsigma_h^m, c_h^m)$ be the solutions from Algorithm 1 and Algorithm 2, respectively. Suppose that $u_0\in L^{\infty}(\Omega; H^2_{per}(D))$ and the conditions of Lemma \ref{Lemma_stability} and Lemma \ref{lemma_fem_stability} satisfy. Then, there exists a constant $C = C(K_1,K_2,K_3,\tilde{K})>0$ such that
		\begin{align}
			\label{error_uh}	&	\max_{1\leq m \leq M}\mE\bigl[\|u^m- u_h^m\|^2_{L^2}\bigr] + \mE\biggl[k\sum_{m=1}^{M}\|\nab(u^{m+\frac12} - u_h^{m+\frac12})\|^2_{L^2}\biggr] \\\nonumber&\qquad\qquad\qquad\qquad\qquad\qquad\leq C\left(k + h^2 + \frac{h^4}{k}\right),\\
			\label{error_sigmah}			& \max_{1\leq m \leq M}\mE\Bigl[\|\vsigma^m - \vsigma_h^m\|^2_{H^1}\Bigr] \leq C\left(k + h^2 + \frac{h^4}{k}\right),\\
			\label{error_ch}	&	\max_{1\leq m \leq M}\mE\bigl[\|c^m- c_h^m\|^2_{L^2}\bigr]  \leq C\left(k + h^2 + \frac{h^4}{k}\right).
		\end{align}
	\end{theorem}	
	\begin{proof}
		First of all, the error estimates \eqref{error_sigmah} and \eqref{error_ch} are a consequence of \eqref{error_uh}. So, it suffices to prove \eqref{error_uh}. 
		
		Denote
		\begin{align*}
			e_{u_h}^m &= u^m - u^m_h  = u^m - \cP_{u} u^m + \cP_{u} u^m - u_h^m:= \theta_{u}^m + \varepsilon_{u}^m  \\
			e_{\vsigma_h}^m &= \vsigma^m - \vsigma_h^m = \vsigma^m - \cP_{\vsigma}\vsigma^m + \cP_{\vsigma} \vsigma^m - \vsigma_h^m:= \vtheta_{\vsigma}^m + \e_{\vsigma}^m. 
		\end{align*}

		Subtracting \eqref{eq_semi_discrete_u} to \eqref{eq_discrete_u}, we obtain the following error equation:
		\begin{align*}
			\bigl(e_{u_h}^{m+1} - e_{u_h}^m,v_h\bigr) + \nu k\bigl(\nab e_{u_h}^{m+\frac12}, \nab v_h\bigr) &= \chi k\bigl(u^{m+\frac12}\vsigma^{m+\frac12} - u_h^{m+\frac12}\vsigma_h^{m+\frac12}, \nab v_h\bigr) \\\nonumber
			&+ \delta\bigl(\vb\cdot\nab e_{u_h}^{m+\frac12}\Delta W_m, v_h\bigr).
		\end{align*}
		
		Taking $v_h = \varepsilon_{u}^{m+\frac12}$, we obtain that
		\begin{align}\label{eq_4.12}
			&\frac12\bigl[\|\varepsilon^{m+1}_u\|^2_{L^2} - \|\varepsilon_u^m\|^2_{L^2}\bigr] + \nu k\|\nab \varepsilon^{m+\frac12}_u\|^2_{L^2}\\\nonumber
			&= -\nu k\bigl(\nab \theta_{u}^{m+\frac12}, \nab \varepsilon_u^{m+\frac12}\bigr) + \chi k\bigl(u^{m+\frac12}e_{\vsigma_h}^{m+\frac12} + e_{u_h}^{m+\frac12}\vsigma_h^{m+\frac12}, \nab \varepsilon_{u}^{m+\frac12}\bigr) \\\nonumber&\qquad + \delta\bigl(\vb\cdot\nab\theta_{u}^{m+\frac12},\varepsilon_{u}^{m+\frac12}\bigr)\Delta W_m  + \delta\bigl(\vb\cdot\nab \varepsilon_{u}^{m+\frac12},\varepsilon_{u}^{m+\frac12}\bigr)\Delta W_{m}\\\nonumber
			&= -\nu k\bigl(\nab \theta_{u}^{m+\frac12}, \nab \varepsilon_u^{m+\frac12}\bigr) + \chi k\bigl(u^{m+\frac12}e_{\vsigma_h}^{m+\frac12} + e_{u_h}^{m+\frac12}\vsigma_h^{m+\frac12}, \nab \varepsilon_{u}^{m+\frac12}\bigr) \\\nonumber&\qquad + \delta\bigl(\vb\cdot\nab\theta_{u}^{m+\frac12},\varepsilon_{u}^{m+\frac12}\bigr)\Delta W_m  \\\nonumber
			&:=I + II + III.
		\end{align}
		
		Using \eqref{ineq_interpolation}, we can control $I$ as below. 
		\begin{align}
			I &\leq \frac{\nu k}{16} \|\nab \varepsilon_{u}^{m+\frac12}\|^2_{L^2} + Ck\|\nab \theta_{u}^{m+\frac12}\|^2_{L^2}  \\\nonumber
			&\leq \frac{\nu k}{16} \|\nab \varepsilon_{u}^{m+\frac12}\|^2_{L^2} + Ckh^2\|\Delta u^{m+\frac12}\|^2_{L^2}.
		\end{align}
		
		To estimate $II$, first we notice that $e_{\vsigma_h}^m$ satisfies the following error equation
		\begin{align*}
			\bigl(e_{\vsigma_h}^{m+1} , \vvarphi_h\bigr) + \bigl(\nab\cdot e_{\vsigma_h}^{m+1},\nab\cdot\vvarphi_h\bigr) + \bigl(rot\, e_{\vsigma_h}^{m+1},\, rot\, \vvarphi_h\bigr) = -\bigl(e_{u}^m, \nab\cdot \vvarphi_h\bigr).
		\end{align*}
		
		With this and the definition of $\cP_{\vsigma}$, we obtain $\forall m \geq 0$
		\begin{align}\label{eq_4.14}
			\bigl(\e_{\vsigma}^{m+1} , \vvarphi_h\bigr) + \bigl(\nab\cdot \e_{\vsigma}^{m+1},\nab\cdot\vvarphi_h\bigr) + \bigl(rot\, \e_{\vsigma}^{m+1},\, rot\, \vvarphi_h\bigr) = -\bigl(e_{u}^m, \nab\cdot \vvarphi_h\bigr),
		\end{align}
		
		which also implies that for all $m\geq 1$
		\begin{align}\label{eq_4.15}
			\bigl(\e_{\vsigma}^{m+\frac12} , \vvarphi_h\bigr) + \bigl(\nab\cdot \e_{\vsigma}^{m+\frac12},\nab\cdot\vvarphi_h\bigr) + \bigl(rot\, \e_{\vsigma}^{m+\frac12},\, rot\, \vvarphi_h\bigr) = -\bigl(e_{u}^{m-\frac12}, \nab\cdot \vvarphi_h\bigr).
		\end{align}
		
		Taking $\vvarphi_h = \e_{\vsigma}^{m+\frac12}$ in \eqref{eq_4.15}, we get
		\begin{align}\label{eq_4.16}
			\|\e_{\vsigma}^{m+\frac12}\|^2_{H^1} &\leq \|e_{u}^{m-\frac12}\|^2_{L^2}\quad\forall m \geq 1.
		\end{align}
		
		Now, we are ready to estimate $II$. Using \eqref{eq_4.16} and \eqref{ineq_interpolation} we have
		\begin{align*}
			II &\leq \frac{\nu k}{16} \|\nab \varepsilon_{u}^{m+\frac12}\|^2_{L^2} + Ck\|u^{m+\frac12}\|^2_{L^4}\|e_{\vsigma_h}^{m+\frac12}\|^2_{L^4} + Ck \|e_{u_h}^{m+\frac12}\|^2_{L^4}\|\vsigma_h^{m+\frac12}\|^2_{L^4}\\\nonumber
			&\leq \frac{\nu k}{16} \|\nab e_{u}^{m+\frac12}\|^2_{L^2} + CK_2 k\|e_{\vsigma_h}^{m+\frac12}\|^2_{L^4} + CK_1k \|e_{u_h}^{m+\frac12}\|^2_{L^4}\\\nonumber
			&\leq \frac{\nu k}{16} \|\nab \varepsilon_{u}^{m+\frac12}\|^2_{L^2} + CK_2 k\|\e_{\vsigma}^{m+\frac12}\|^2_{L^4} +CK_2 k\|\vtheta_{\vsigma}^{m+\frac12}\|^2_{L^4} \\\nonumber
			&\qquad + CK_1k \|\varepsilon_{u}^{m+\frac12}\|^2_{L^4}  + CK_1k \|\theta_{u}^{m+\frac12}\|^2_{L^4}\\\nonumber
			&\leq \frac{\nu k}{16} \|\nab \varepsilon_{u}^{m+\frac12}\|^2_{L^2} + CK_2 k\|\e_{\vsigma}^{m+\frac12}\|^2_{H^1} +CK_2 kh^3\|{\vsigma}^{m+\frac12}\|^2_{H^2} \\\nonumber
			&\qquad + CK_1k \|\varepsilon_{u}^{m+\frac12}\|_{L^2}\|\nab \varepsilon_{u}^{m+\frac12}\|_{L^2}  + CK_1kh^3 \|{u}^{m+\frac12}\|^2_{H^2}\\\nonumber
			&\leq \frac{\nu k}{16} \|\nab \varepsilon_{u}^{m+\frac12}\|^2_{L^2} + \frac{\nu k}{16} \|\nab \varepsilon_{u}^{m+\frac12}\|^2_{L^2} + CK_2 k\|e_{u}^{m-\frac12}\|^2_{L^2} +CK_2 kh^3\|{\vsigma}^{m+\frac12}\|^2_{H^2} \\\nonumber
			&\qquad + CK^2_1k \|\varepsilon_{u}^{m+\frac12}\|^2_{L^2}  + CK_1kh^3 \|{u}^{m+\frac12}\|^2_{H^2}\\\nonumber
			&\leq \frac{\nu k}{16} \|\nab \varepsilon_{u}^{m+\frac12}\|^2_{L^2} + \frac{\nu k}{16} \|\nab \varepsilon_{u}^{m+\frac12}\|^2_{L^2} + CK_1^2 k\|\varepsilon_{u}^{m+\frac12}\|^2_{L^2} + CK_2 k\|\varepsilon_{u}^{m-\frac12}\|^2_{L^2} \\\nonumber
			&\qquad+CK_2 k h^4\| u^{m+\frac12}\|^2_{H^2}  +  CK_2 kh^3\|{\vsigma}^{m+\frac12}\|^2_{H^2}  + CK_1kh^3 \|{u}^{m+\frac12}\|^2_{H^2}.
		\end{align*}
		
		Next, we estimate the noise error terms. Using integration by parts, and \eqref{ineq_interpolation}, we obtain
		\begin{align*}
			III &= -\delta\bigl(\theta_{u}^{m+\frac12},\vb\cdot\nab \varepsilon_{u}^{m+\frac12}\bigr)\Delta W_m\\\nonumber
			&\leq \delta |\vb|\|\theta_{u}^{m+\frac12}\|_{L^2}\|\nab\varepsilon_{u}^{m+\frac12}\|_{L^2}|\Delta W_m|\\\nonumber
			&\leq \frac{\nu k}{16} \|\nab \varepsilon_{u}^{m+\frac12}\|^2_{L^2} + \frac{\delta^2|\vb|^2}{\nu} \|\theta_{u}^{m+\frac12}\|^2_{L^2}\frac{|\Delta W_m|^2}{k}\\\nonumber
			&\leq \frac{\nu k}{16} \|\nab \varepsilon_{u}^{m+\frac12}\|^2_{L^2} + \frac{C\delta^2|\vb|^2}{\nu} h^4\|{u}^{m+\frac12}\|^2_{H^2}\frac{|\Delta W_m|^2}{k},
		\end{align*}
		which together with Lemma \ref{Lemma_stability} implies that
		\begin{align}
			\mE[III] \leq  \frac{\nu k}{16}\mE[\|\nab \varepsilon_{u}^{m+\frac12}\|^2_{L^2}] + \frac{C\delta^2|\vb|^2K_3}{\nu} h^4
		\end{align}
		
		Next, applying the summation $\sum_{m=1}^{\ell}$ for $1\leq \ell\leq M-1$,  taking the expectation and then substituting all the estimations of $I, II ,III$ to \eqref{eq_4.12}, we obtain
		\begin{align*}
			&	\frac12\mE[\|\varepsilon_{u}^{\ell +1}\|^2_{L^2}] + \frac{3\nu k}{4}\sum_{m=1}^{\ell} \mE[\|\nab\varepsilon_{u}^{m+\frac12}\|^2_{L^2}]\\\nonumber
			&\leq Ch^2(1 + K_1 + K_2)k\sum_{m=1}^{\ell} \mE[\| u^{m+\frac12}\|^2_{H^2}]  + Ch^2K_2k\sum_{m=1}^{\ell}\mE\bigl[\|\vsigma^{m+\frac12}\|^2_{H^2}\bigr]\\\nonumber
			&\qquad+ \frac{C\delta^2 |\vb|^2K_3}{\nu}\frac{h^4}{k} + C(K_1^2 + K_2) k\sum_{m=1}^{\ell}\mE[\|\varepsilon_{u}^{m-\frac12}\|^2_{L^2}] + \frac12\mE[\|\varepsilon_{u}^1\|^2_{L^2}]\\\nonumber
			&\leq Ch^2(1 + K_1 + K_2)K_2 + Ch^2K_2^2+ \frac{C\delta^2 |\vb|^2K_3}{\nu}\frac{h^4}{k} \\\nonumber
			&\qquad+ C(K_1^2 + K_2) k\sum_{m=1}^{\ell}\mE[\|\varepsilon_{u}^{m-\frac12}\|^2_{L^2}] + \frac12\mE[\|\varepsilon_{u}^1\|^2_{L^2}].
		\end{align*}
		
		Using the standard discrete Gronwall inequality, we obtain
		\begin{align}\label{eq_4.17}
			&	\frac12\mE[\|\varepsilon_{u}^{\ell +1}\|^2_{L^2}] + \frac{3\nu k}{4}\sum_{m=1}^{\ell} \mE[\|\nab\varepsilon_{u}^{m+\frac12}\|^2_{L^2}]\\\nonumber
			&\leq \Bigl[C(K_2+K_1K_2+  K^2)h^2 + \frac{CK_3|\vb|^2\delta^2}{\nu}\frac{h^4}{k}\\\nonumber
			&\qquad\qquad\qquad+\frac12\mE[\|\varepsilon_{u}^1\|^2_{L^2}]\Bigr] \exp(C(K^2_1 + K_2)T).
		\end{align}
		
		Finally, it is left to estimate $\varepsilon_{u}^1$. Taking $m=0$ in \eqref{eq_4.12}, then using the assumptions that $\varepsilon_{u}^0=0, \e_{\vsigma}^0=0$, we get that
		\begin{align}\label{eq_4.19}
			\|\varepsilon_{u}^1\|^2_{L^2} +\frac{\nu k}{2} \|\nab \varepsilon_{u}^{1}\|^2_{L^2} &=  -\frac{ \nu k}{2} \bigl(\nab \theta_{u}^{1},\nab \varepsilon_{u}^{1}\bigr)   + \frac12\chi k\bigl(u^{\frac12} e_{\vsigma_h}^{1} + e_{u_h}^1\vsigma_h^{\frac12},\nab \varepsilon_{u}^{1}\bigr) \\\nonumber&\qquad+  \frac{\delta}{2}\bigl(\vb\cdot\nab\theta_{u}^{1},\varepsilon_{u}^1\bigr)\Delta W_0\\\nonumber
			&\leq \frac{\nu k}{16} \|\nab \varepsilon_{u}^1\|^2_{L^2} + Ckh^2\|{u}^{1}\|^2_{H^2}\\\nonumber
			&\qquad + Ck h^3 \|u^{\frac12}\|^2_{H^1}\|\vsigma^1\|^2_{H^2} + Ckh^4\|u^{\frac12}\|^2_{H^1}\|u^0\|^2_{H^2}\\\nonumber
			&\qquad+ Ckh^3 \|u^1\|^2_{H^2}\|\vsigma_h^{\frac12}\|^2_{H^1} + Ck\|\vsigma_h^{\frac12}\|^2_{H^1}\|\varepsilon_{u}^1\|^2_{H^1}\\\nonumber
			&\qquad + \frac14\|\varepsilon_{u}^1\|^2_{L^2} + C\delta^2|\vb|^2 h^2\|u^1\|^2_{H^2}|\Delta W_0|^2\\\nonumber
			&\leq \frac{\nu k}{16} \|\nab \varepsilon_{u}^1\|^2_{L^2} + CK_3kh^2 + \frac14\|\varepsilon_{u}^1\|^2_{L^2} \\\nonumber
			&\qquad + CK^2_2h^3 + CK_1h^4\|u^0\|^2_{H^2}\\\nonumber
			&\qquad+ CK_3K_1kh^3 + CK_1(K_1 + \tilde{K})k
			+ CK_3 h^2\\\nonumber
			&=\frac{\nu k}{16} \|\nab \varepsilon_{u}^1\|^2_{L^2}  + \frac14\|\varepsilon_{u}^1\|^2_{L^2} +  CK_1(K_1 + \tilde{K})k\\\nonumber
			&\qquad +  C(kK_3 + hK_2^2 + h^2K_1\|u_0\|^2_{H^2}+ k K_1K_3 + K_3)h^2,
		\end{align}
		which implies that
		\begin{align}\label{eq_4.21}
			\mE[\|\varepsilon_{u}^1\|^2_{L^2}] + \nu k\mE[\|\nab \varepsilon_u^1\|^2_{L^2}] \leq C(k+h^2).
		\end{align}
		
		The proof is complete by combining \eqref{eq_4.17} and \eqref{eq_4.21}.
		
	\end{proof}

	Finally, we present the global error estimates in the following theorem. The proof of this theorem can be derived directly by the triangle inequality and Theorem \ref{Thm_error_semi}, and Theorem \ref{Thm_error}.
	\begin{theorem}\label{Thm_global_error} Let $(u_h^m, \vsigma_h^m, c_h^m)$ be the solution from Algorithm 2. Suppose that $u_0\in L^{\infty}(\Omega; H^2_{per}(D))$ and the conditions of Theorem \ref{Thm_error_semi} and Theorem \ref{Thm_error} satisfy. Then, there exists a constant $C>0$ such that
		\begin{align*}
			\max_{1\leq m \leq M}\mE\bigl[\|u(t_m)- u_h^m\|^2_{L^2}\bigr]+ \mE\biggl[k\sum_{m=1}^{M}\|\nab(u(t_{m+\frac12}) - u_h^{m+\frac12})\|^2_{L^2}\biggr]
			&\leq C\left(k + h^2 + \frac{h^4}{k}\right),\\
			\max_{1\leq m \leq M}	\mE\Bigl[\|\vsigma(t_m) - \vsigma_h^m\|^2_{H^1}\Bigr] &\leq C\left(k + h^2 + \frac{h^4}{k}\right),\\
			\max_{1\leq m \leq M}\mE\bigl[\|c(t_m) - c_h^m\|^2_{L^2}\bigr]  &\leq C\left(k + h^2 + \frac{h^4}{k}\right).
		\end{align*}
	\end{theorem}

	\section{Numerical experiments}\label{sec5}
	In this part, we present a series of numerical tests to verify the error bounds of Theorem \ref{Thm_global_error}. In all our experiments   we set $D = (0,1)^2\subset \mathbb{R}^2$, $T =1, \nu =1$ and $\chi =1$.
	We choose $W(t)$ in \eqref{eq1.1} to be a Wiener process with $\mathbb{R}$ value that is simulated by the minimal time-step size $k_0 = 1/2048$. For all tests, we use the standard Monte Carlo method to compute the expectation. Moreover, in Tests 1, 2, and 3, we choose $u_0 = \sin(\pi x)\sin(\pi y)$ and $c_0 = \sin(\pi x)\sin(\pi y)$.

	We implement Algorithm 2 and compute the errors of the approximate solutions $(\vsigma^m_h,u^m_h,c_h^m)$ in the norms specified below. Since exact solutions are unknown, errors are calculated between the computed solution $(\vsigma^m_h(\omega_j), u_h^m(\omega_j), c_h^m(\omega_j))$ and a reference solution $(\vsigma^m_{ref}(\omega_j),u^m_{ref}(\omega_j), {{c}}^m_{ref}(\omega_j))$ (specified later) in the $\omega_j$-th sample.
	
	Furthermore, to evaluate errors in strong norms, we use the following numerical integration formulas: 
	\begin{align*}
		L^2_{\omega}L^{\infty}_tL^2_x(u)&:=	\bigl(\mE\bigl[\max_{1\leq m \leq M}\|u(t_m) - u^m_h\|^2_{L^2}\bigr]\bigr)^{1/2} \\
		&\approx \Bigl(\frac{1}{J}\sum_{j= 1}^J\bigl(\max_{1\leq m \leq M} \|u^m_{ref}(\omega_j) - u^m_h(\omega_j)\|^2_{L^2}\bigr)\Bigr)^{1/2},\\
		L^2_{\omega}L^{\infty}_tL^2_x(c)&:=	\bigl(\mE\bigl[\max_{1\leq m \leq M}\|c(t_m) - c^m_h\|^2_{L^2}\bigr]\bigr)^{1/2} \\
		&\approx \Bigl(\frac{1}{J}\sum_{j= 1}^J\bigl(\max_{1\leq m \leq M} \|c^m_{ref}(\omega_j) - c^m_h(\omega_j)\|^2_{L^2}\bigr)\Bigr)^{1/2},\\
		L^2_{\omega}L^{\infty}_tH^1_x(\vsigma)&:=	\bigl(\mE\bigl[\max_{1\leq m \leq M}\|\vsigma(t_m) - \vsigma^m_h\|^2_{H^1}\bigr]\bigr)^{1/2} \\
		&\approx \Bigl(\frac{1}{J}\sum_{j= 1}^J\bigl(\max_{1\leq m \leq M} \|\vsigma^m_{ref}(\omega_j) - \vsigma^m_h(\omega_j)\|^2_{L^2}\bigr)\Bigr)^{1/2},
	\end{align*}
	where $J>0$ is the number of samples from the Monte Carlo method. 
	
	{\bf Test 1.} In this test, we aim to verify the global convergence rates of the approximate solutions \( (\vsigma_h^m, u_h^m, c_h^m) \) as predicted by Theorem~\ref{Thm_global_error}. The theoretical error bound \( \mathcal{O}(k^{1/2} + h + k^{-1/2}h^2) \) contains the term \( k^{-1/2} \), which can dominate if the spatial and temporal discretizations are not properly balanced. To achieve a global convergence rate of \( \mathcal{O}(h) \), we set the time step size to \( k = h^2 \). 
	
	To validate this convergence rate, we run Algorithm~2 with mesh sizes \( h = 1/2, 1/4, 1/8, 1/16 \) and corresponding time steps \( k = 1/4, 1/16, 1/64, 1/256 \). The reference solutions \( \{ (\vsigma_{\mathrm{ref}}^m, u_{\mathrm{ref}}^m, c_{\mathrm{ref}}^m) \} \) are computed using finer discretizations \( (h_{\mathrm{ref}}, k_{\mathrm{ref}}) = (h/2, k/4) \). That is, the errors are approximated by comparing numerical solutions on two consecutive spatial grids \cite{feng2020fully, vo2023higher, FengVo2024, vo2025high}.
	
\begin{figure}[htp]
	\begin{center}
			\includegraphics[scale=0.11]{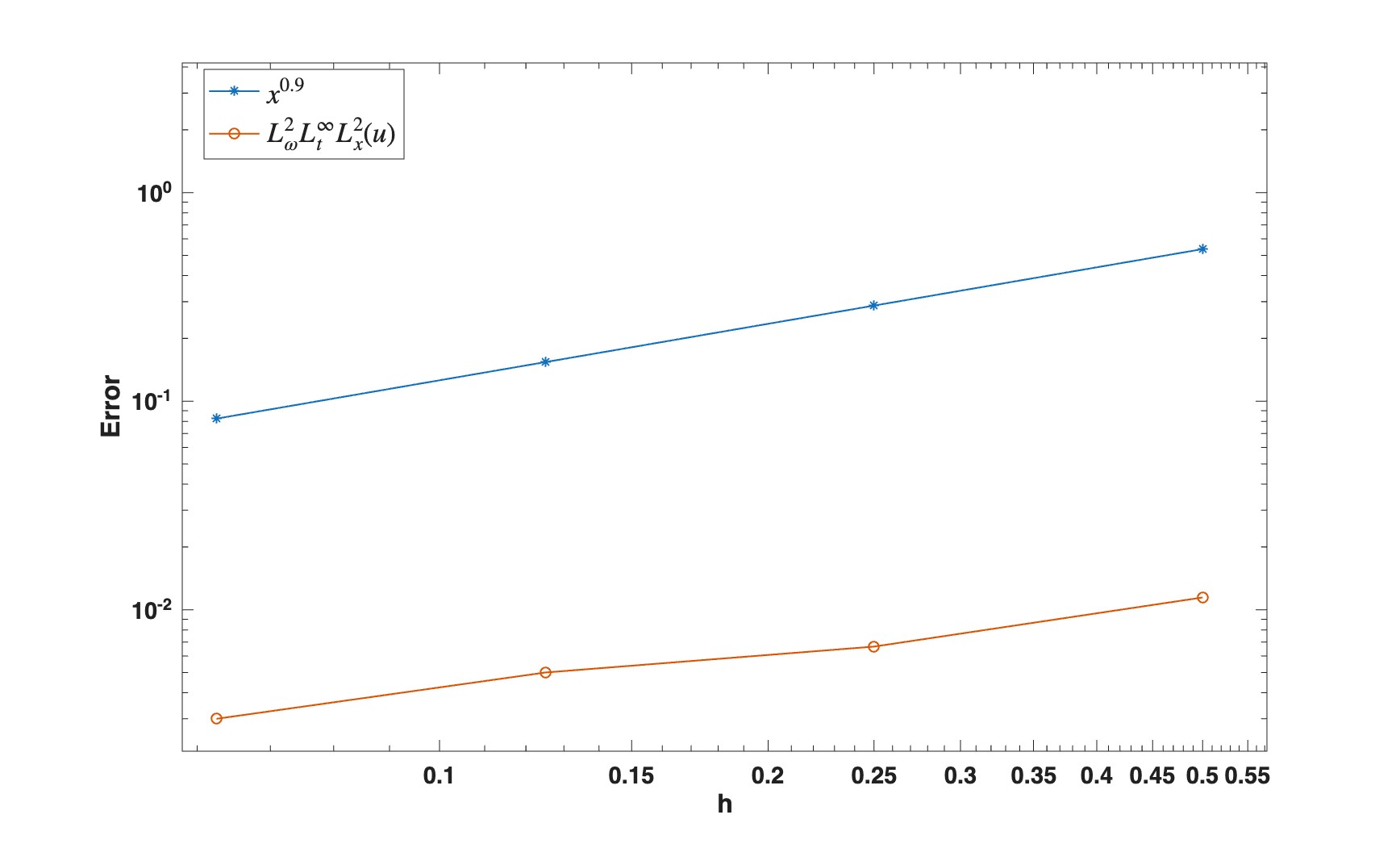}
		\includegraphics[scale=0.11]{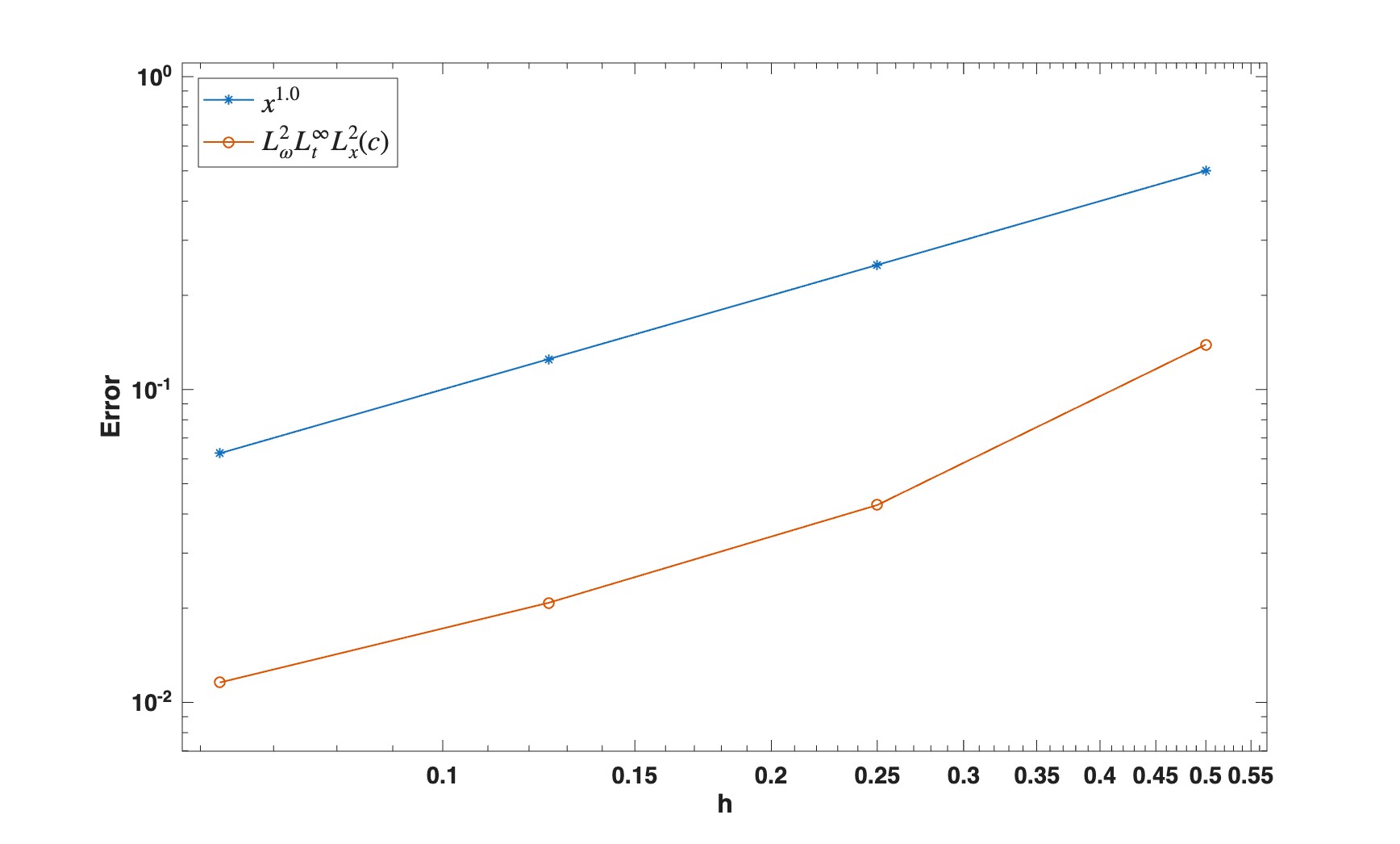}
		\includegraphics[scale=0.11]{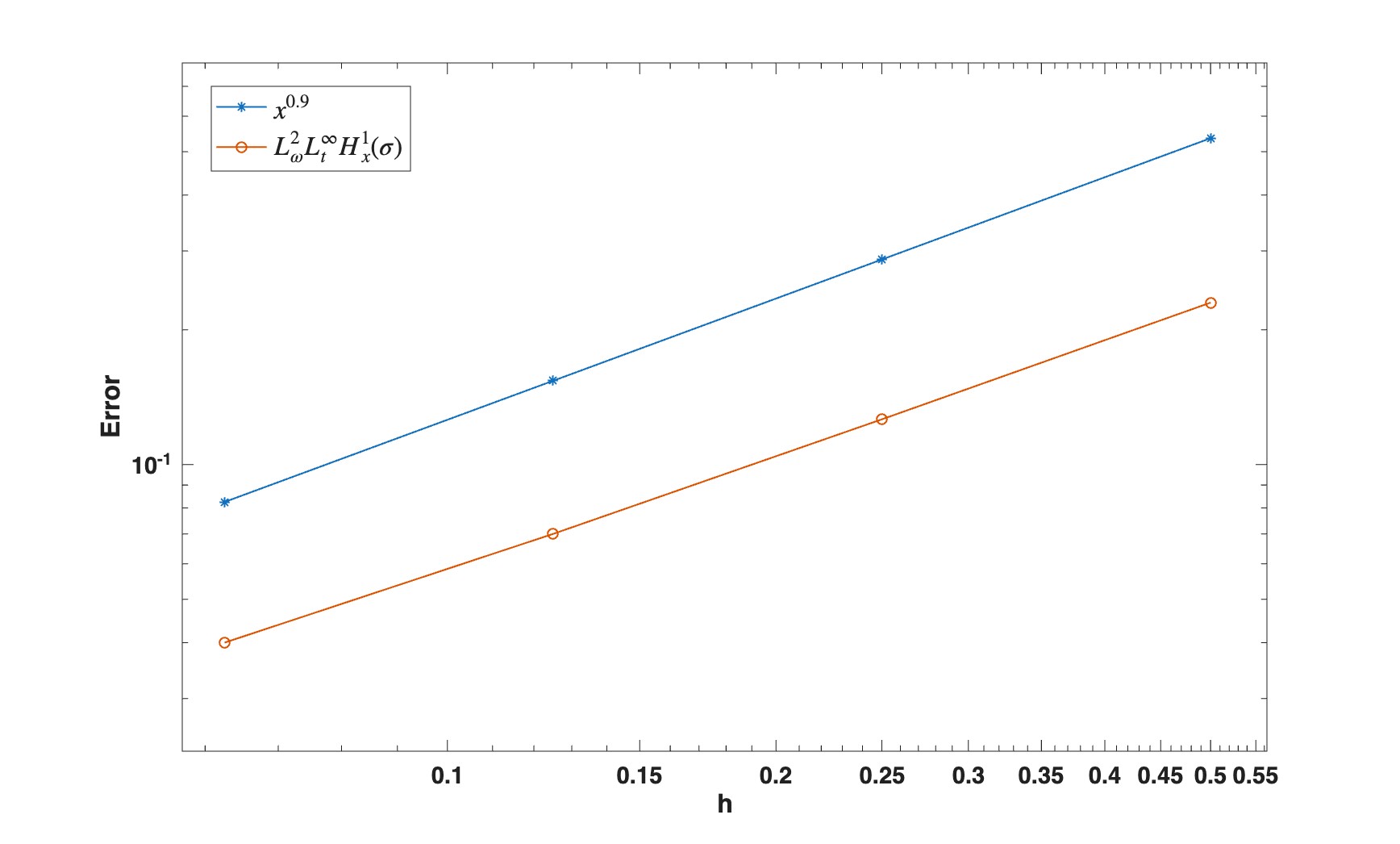}
		\caption{Plots of the errors and convergence order of the computed solution: Top left: $u_h^m$; Top right: $c_h^m$; Bottom: $\vsigma_h^m$, with $k = h^2$, $\delta = 1$, $J = 400$.}\label{fig5.1}
			\end{center}
	\end{figure}

	The computed errors and the observed convergence rates are illustrated in Figure~\ref{fig5.1}. As expected, the numerical results confirm a first-order convergence rate for all components, in agreement with Theorem~\ref{Thm_global_error}. These results also highlight the influence of noise on convergence behavior, leading to slower rates compared to the deterministic setting.
	
	\medskip

\begin{figure}[htp]
		\begin{center}
			\includegraphics[scale=0.12]{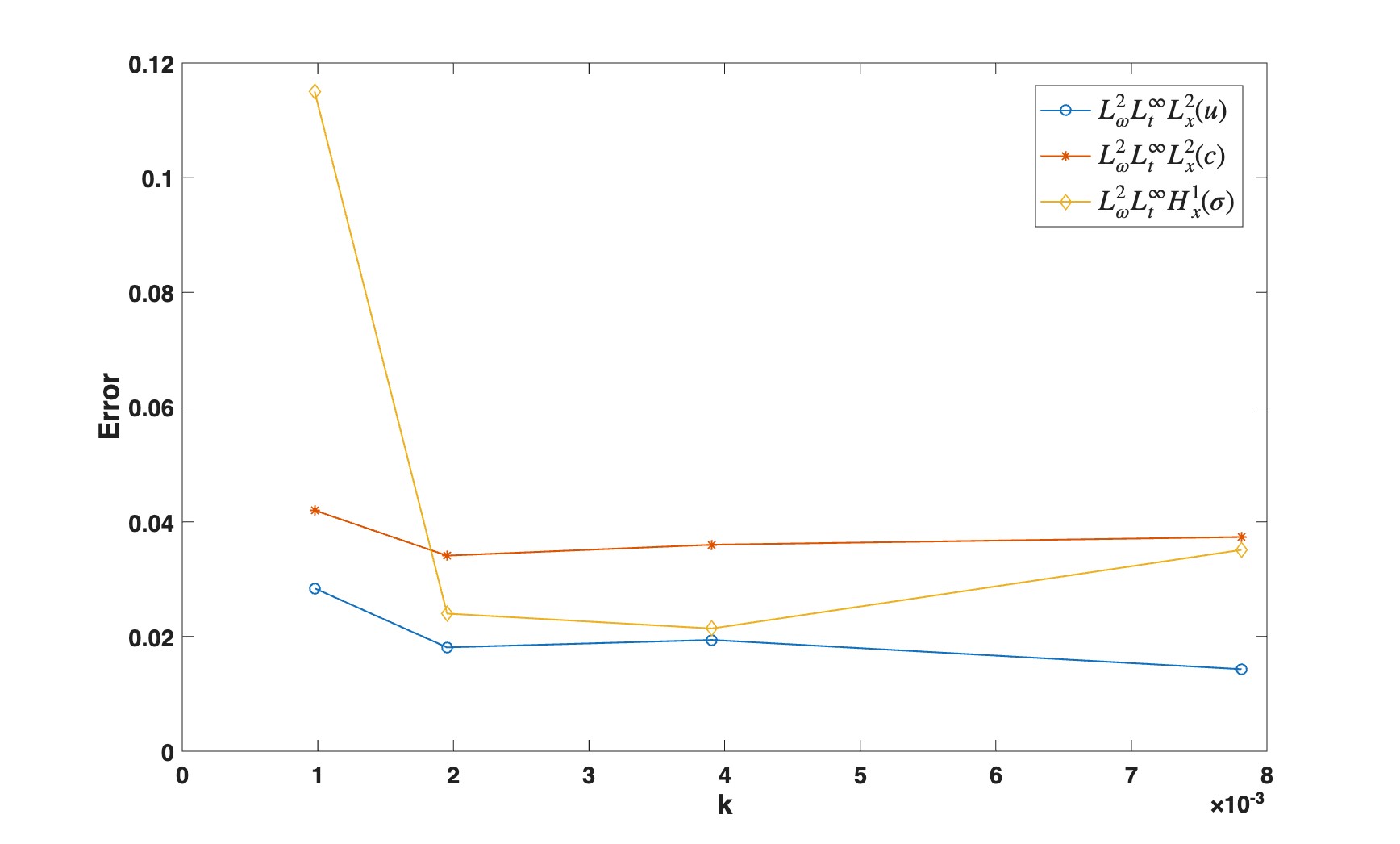}
			\caption{{Plots of the errors of the computed solutions $\{\vsigma^m_h, { u}^m_h, c^m_h\}$ for $k \ll h^2$, and $\delta = 10$, $J = 400$.}}\label{fig5.2}
		\end{center}
	\end{figure}

	{\bf Test 2.} In this test, we aim to investigate the emergence of the factor \( k^{-\frac{1}{2}} \) in the error bounds. To this end, we run Algorithm~2 using a fixed mesh size \( h = 1/10 \) while successively refining the time step size with \( k = 1/128, 1/256, 1/512, 1/1024 \). Unlike in Test~1, the relationship \( k = h^2 \), which ensures a balance between spatial and temporal discretization, is violated here. That means we chose $k \ll h^2$. As a result, we expect the errors to grow as the time step size decreases. The errors of the calculated solutions \( \{ \vsigma_h^m, u_h^m, c_h^m \} \) are presented in Figure~\ref{fig5.2}. As anticipated, the numerical results show that the errors increase as \( k \) becomes smaller, thereby confirming the predicted influence of the factor \( k^{-\frac{1}{2}} \) in the error bounds.
	
\begin{figure}[htp]
	\begin{center}
			\includegraphics[scale=0.1]{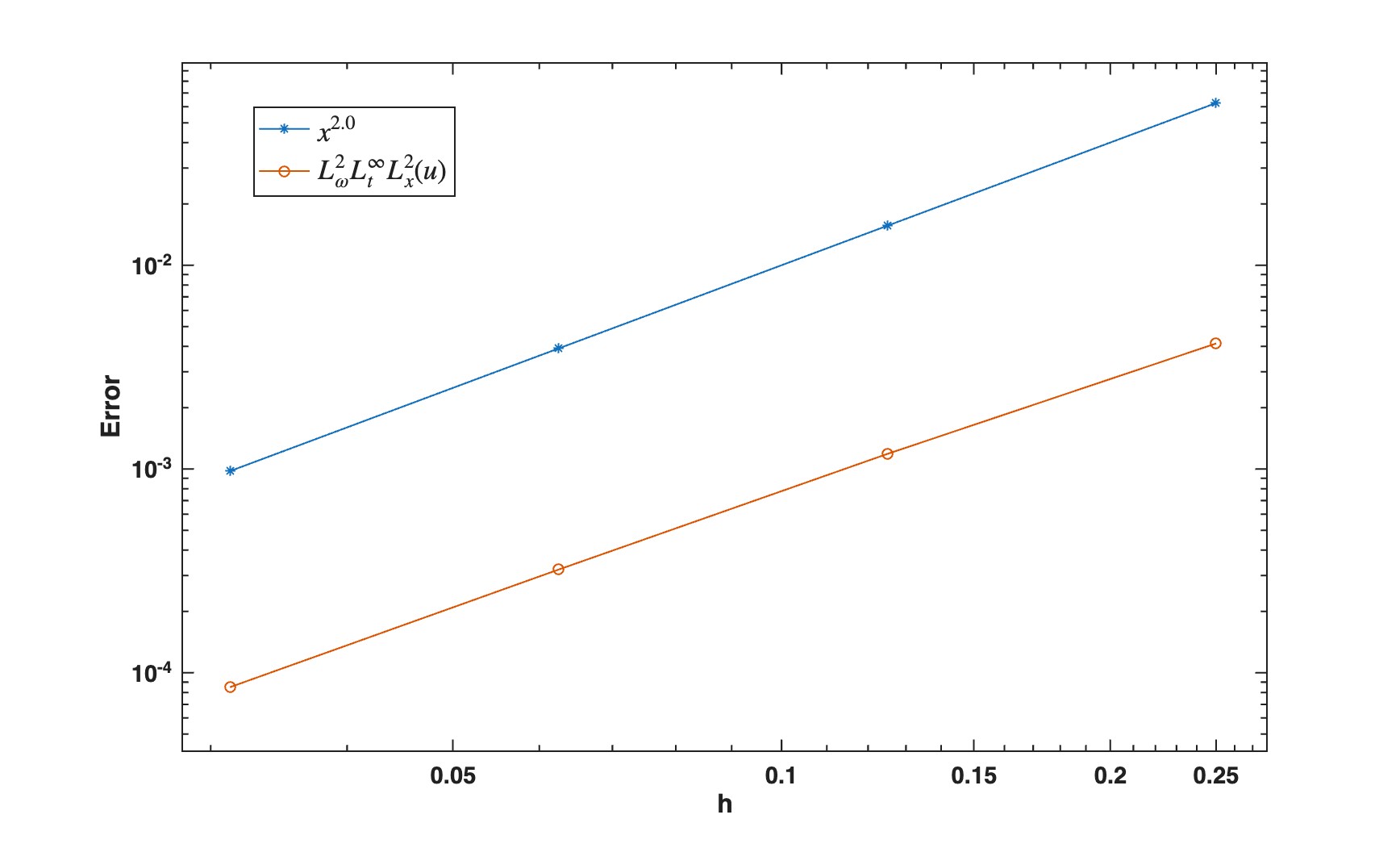}
		\includegraphics[scale=0.1]{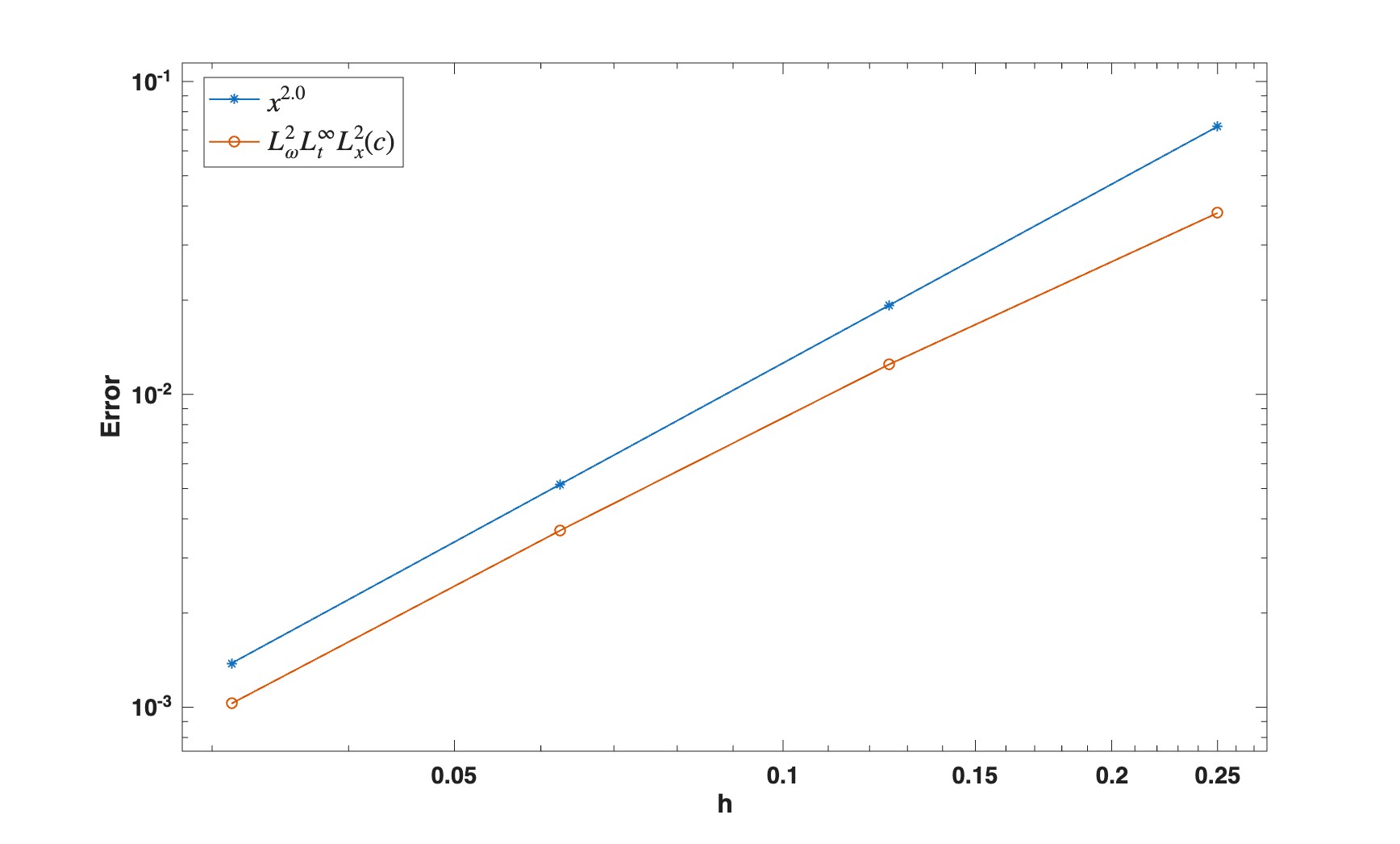}
		\includegraphics[scale=0.1]{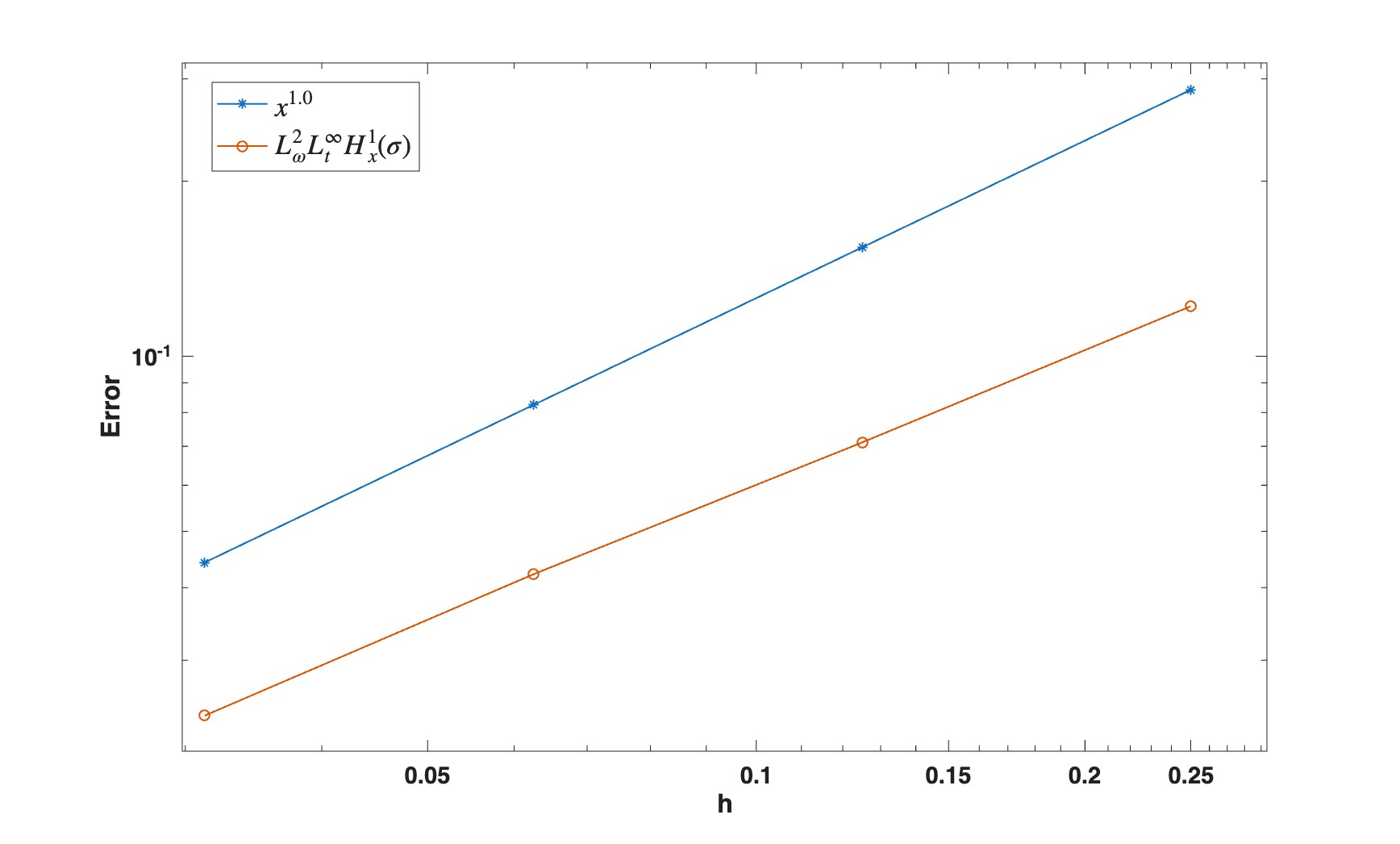}
		\caption{Plots of the errors and convergence order of the computed solution: Top left: $u_h^m$; Top right: $c^m_h$; Bottom: $\vsigma_h^m$, with $\delta = 0.1$, $J = 400$.}\label{fig5.3}
	\end{center}
	\end{figure}

	\medskip
	{\bf Test 3.} In this experiment, we investigate how the noise intensity influences the convergence rates of Algorithm~2. Specifically, we consider a small noise level of \( \delta = 0.1 \), which is significantly smaller than those used in Tests~1 and~2. At this level, the system~\eqref{eq1.1} behaves almost deterministically, and we therefore expect Algorithm~2 to exhibit convergence rates comparable to those in the deterministic setting. We fix the time step \( k = 1/2048 \) and vary the spatial mesh size \( h = 1/4, 1/8, 1/16, 1/32 \) to compute the errors of the approximations \( \{ \vsigma_h^m, u_h^m, c_h^m \} \). The convergence behavior is illustrated in Figure~\ref{fig5.3}, which shows that \( u_h^m \) and \( c_h^m \) achieve second-order accuracy in the \( L^2 \)-norm, while \( \vsigma_h^m \) converges with first-order accuracy in the \( H^1 \)-norm. These results are consistent with the deterministic case and highlight the impact of noise intensity on the convergence properties of Algorithm~2.
	
	\medskip
	{\bf Test 4.} In this numerical experiment, we aim to verify two key properties of the stochastic Keller--Segel system: the \emph{finite-time blow-up} of the solution and the \emph{positivity} of the approximated density \( u_h^m \). We choose parameters \( \chi = 4\pi \), \( \nu = 1 \), and \( \delta = 1.0 \), with initial conditions
	\[
	u_0(x, y) = 1000\, e^{-100(x^2 + y^2)}, \quad c_0(x, y) = 500\, e^{-50(x^2 + y^2)},
	\]
	on the domain $D = [-\frac12,\frac12]\times[-\frac12,\frac12]$.
	The expectation \( \mathbb{E}[u_h^M] \) is computed using a spatial mesh size \( h = \frac{1}{60} \) and a time step \( k = 10^{-6} \).
	
	Our numerical results indicate that the solution is positive and exhibits blow-up behavior as the final time increases, with \( t_M = 3 \times 10^{-5},\ 5 \times 10^{-5},\ 9 \times 10^{-5},\ 2 \times 10^{-4} \). Additionally, the presence of stochastic noise appears to \emph{accelerate the blow-up}, leading to qualitative differences compared to the deterministic setting. This highlights the significant role of noise in the dynamics of chemotactic aggregation.

	\begin{figure}[htp]
		\includegraphics[scale=0.12]{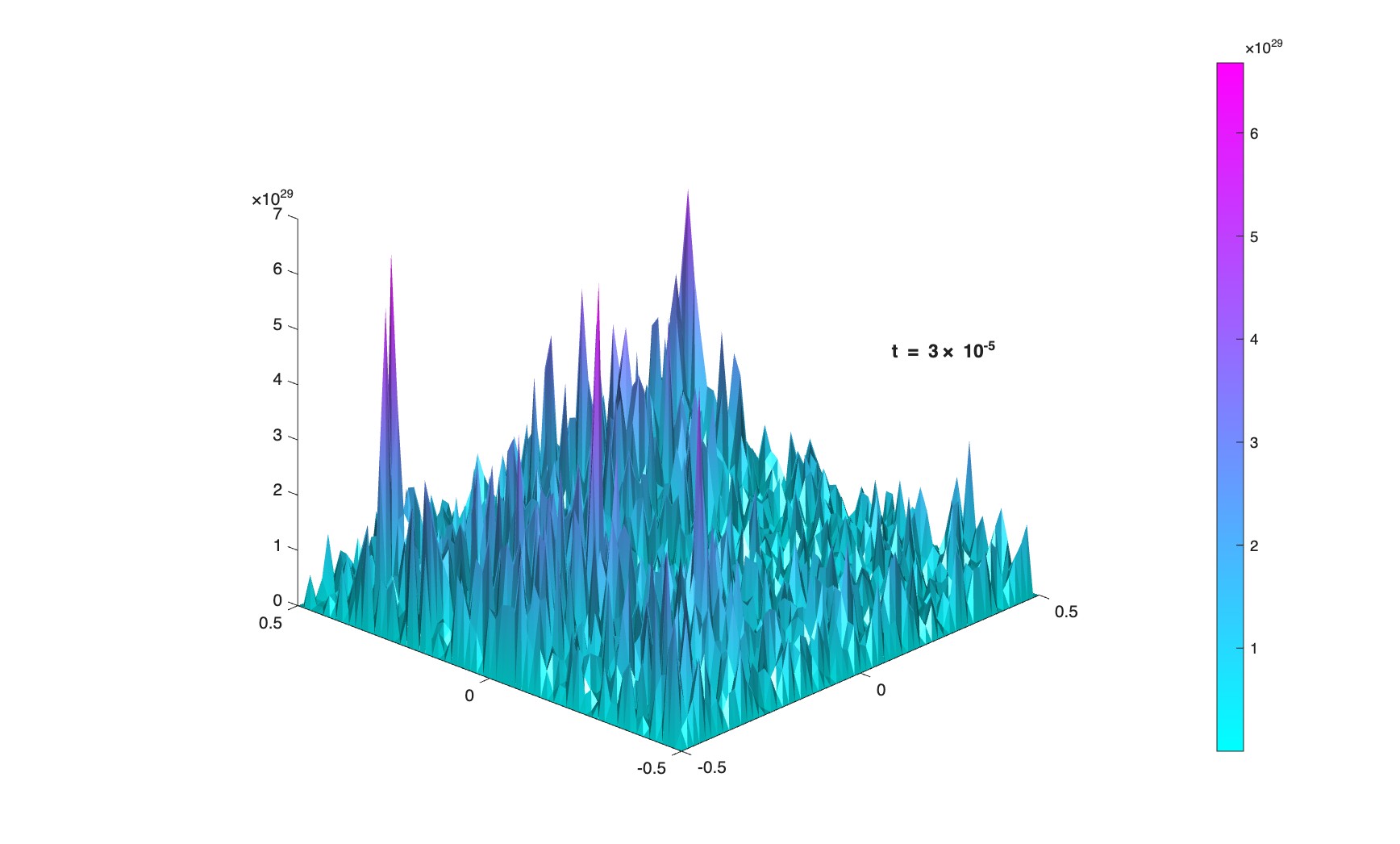}
		\includegraphics[scale=0.12]{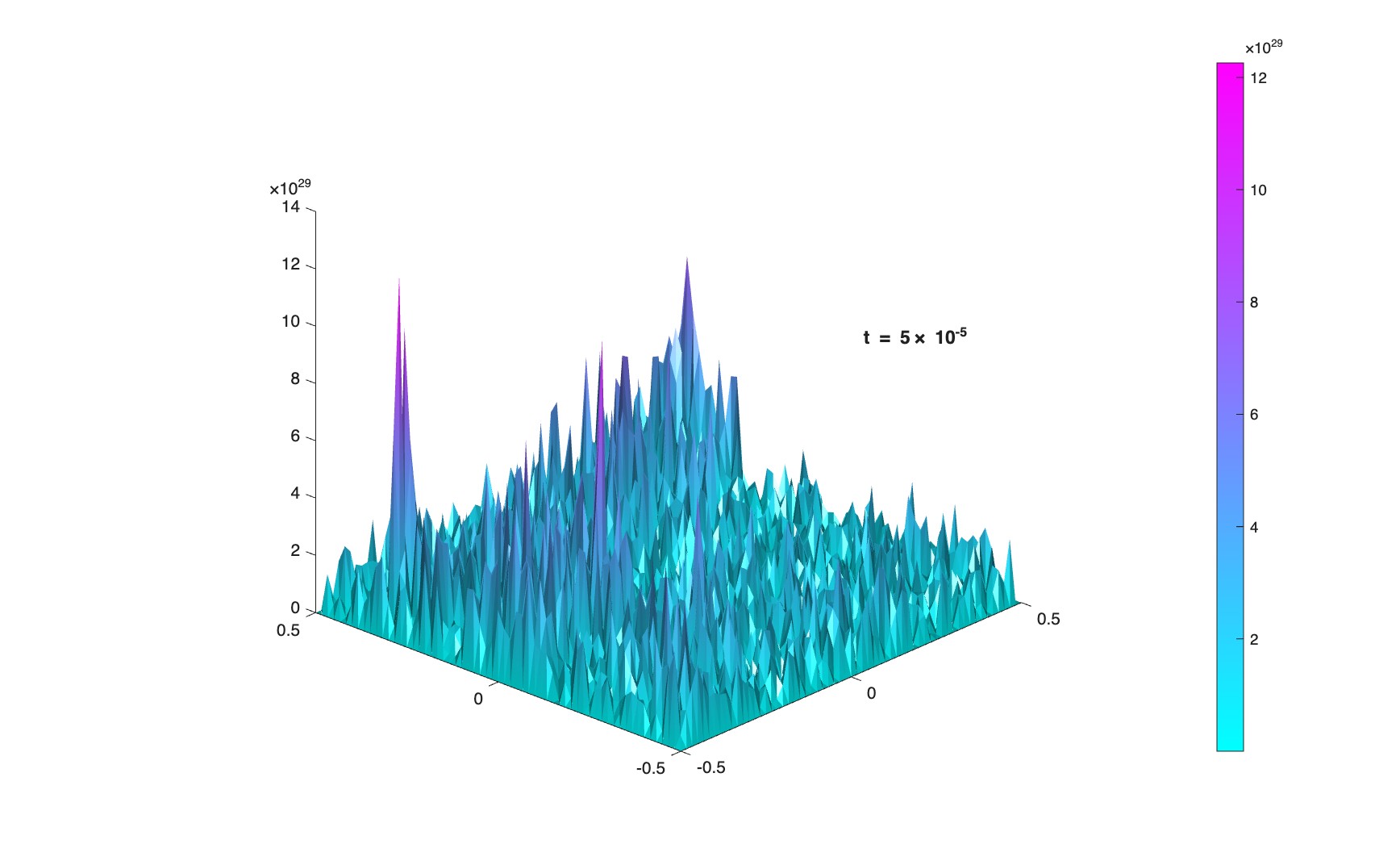}
		\includegraphics[scale=0.12]{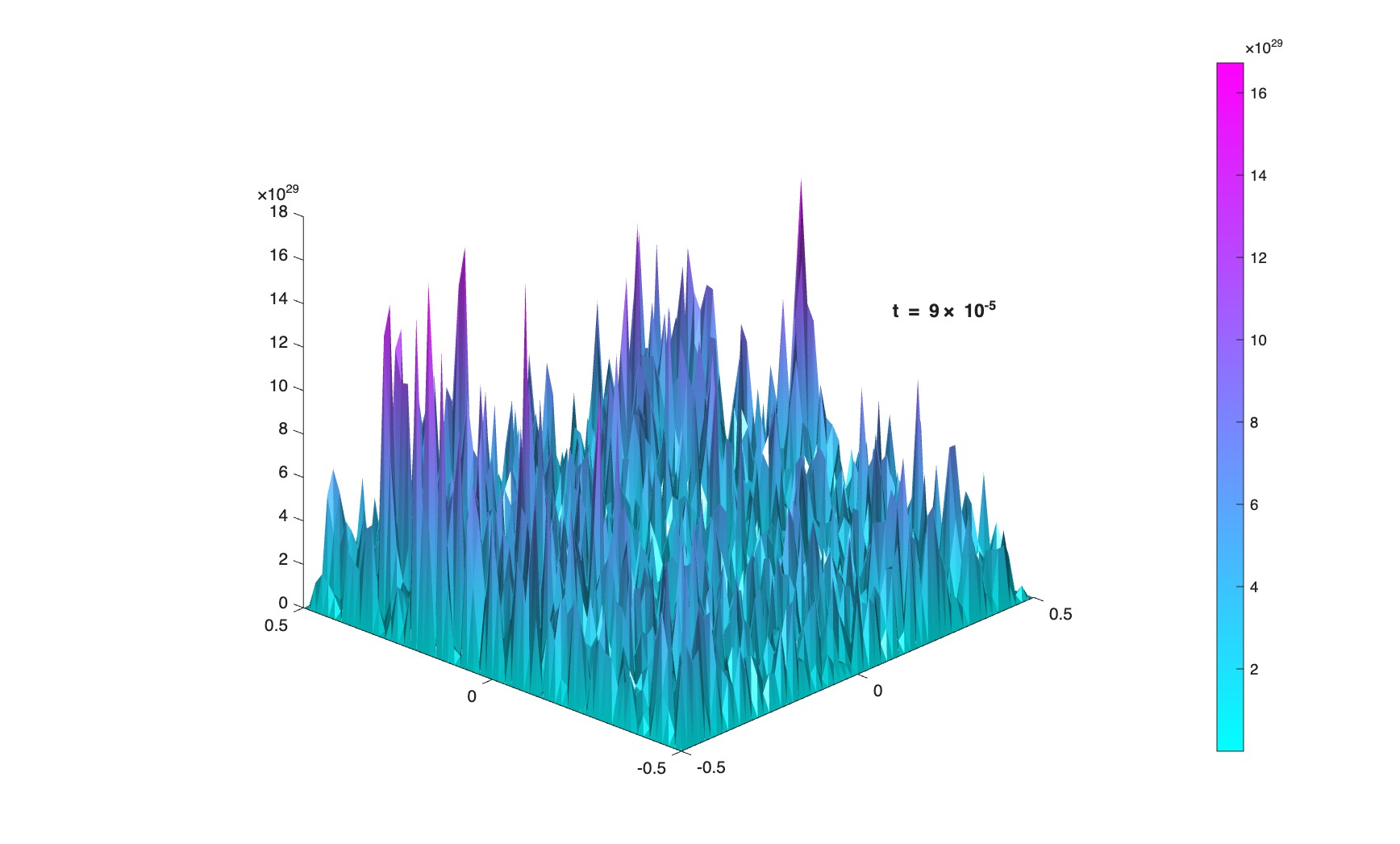}
		\includegraphics[scale=0.12]{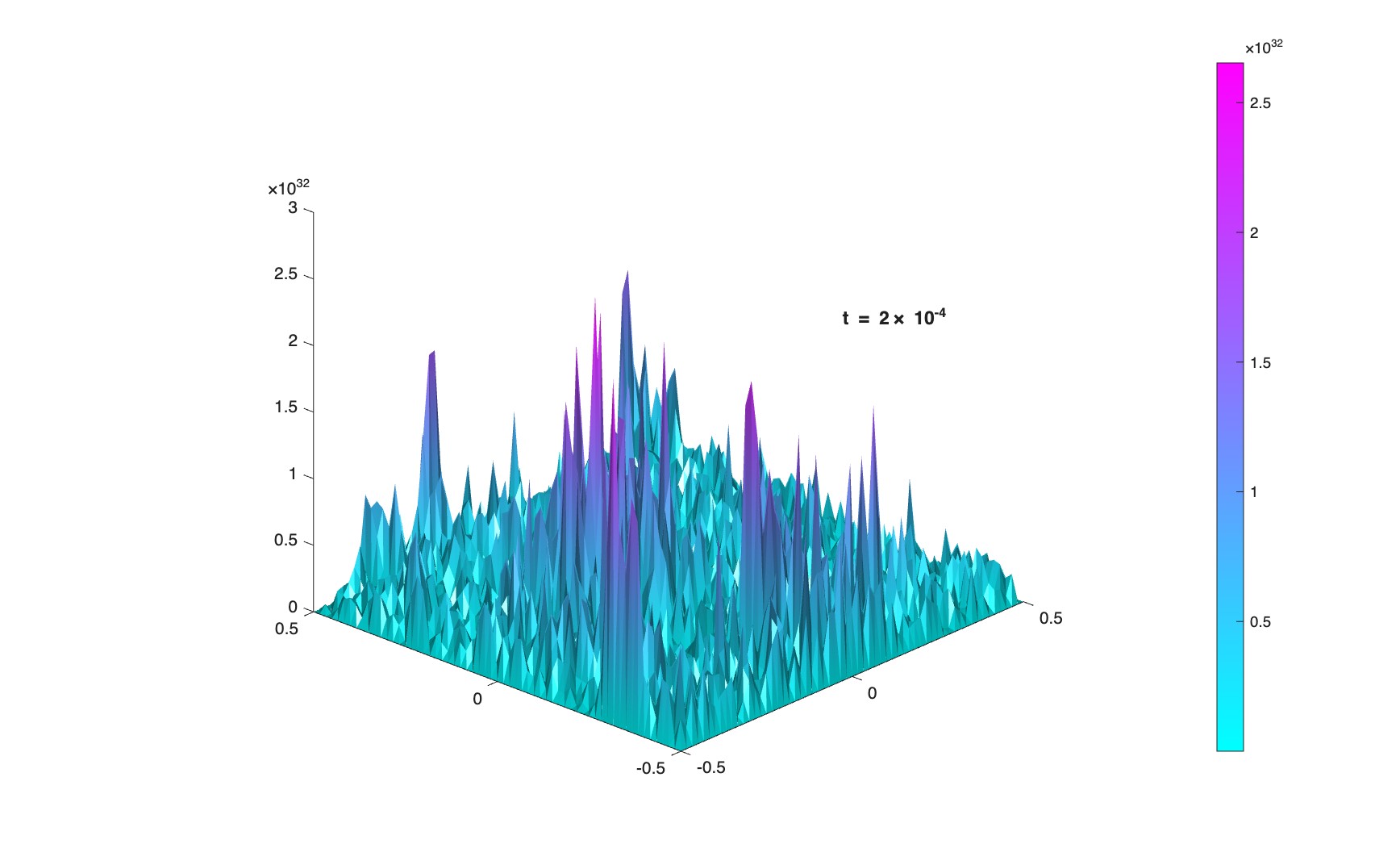}
		\caption{Blowup plots of the computed solution $\mE[u_h^M]$ with different final times $t_M = 3\times 10^{-5}, 5\times 10^{-5}, 9\times 10^{-5}, 2\times 10^{-4}$, and $J = 400$.}\label{fig5.5}
	\end{figure}
	
	\section{Conclusion}\label{sec6} In this paper, we presented numerical methods for solving the stochastic Keller–Segel system. This work serves as a foundational result and provides a reference framework for developing numerical schemes with rigorous error analysis for such stochastic PDEs. A key finding is that the presence of stochastic noise fundamentally alters the behavior of the error bounds, making them depend inversely on the time step size, which is a significant departure from the deterministic case.
	
	All error estimates in this study are derived under the assumption that the initial data are deterministic. Extending these results to the case of initial random data presents additional challenges. In future work, we plan to address this by employing local sample space techniques \cite{carelli2012rates} to manage nonlinearities. This strategy is expected to yield convergence results in probability, offering a weaker convergence compared with the current results in the paper.

	\bibliographystyle{abbrv}
	\bibliography{references}

\end{document}